\documentclass[oneside,final,reqno]{amsart}
\usepackage{geometry}

\usepackage{graphicx}
\usepackage{amsmath}
\usepackage[utf8]{inputenc}
\usepackage{nicefrac}
\usepackage{mathabx}
\usepackage{subcaption}
\usepackage{amsrefs}
\usepackage{esint}
\usepackage{mathtools}
\usepackage{graphicx}

\usepackage{textcomp} 

\usepackage{wrapfig}
\usepackage{xcolor}
\usepackage{enumitem}
\usepackage{verbatim}
\usepackage{hyperref}
\usepackage[right]{showlabels}
\usepackage{cleveref}
\expandafter\def\csname ver@etex.sty\endcsname{3000/12/31}

\usepackage{autonum}
\usepackage[normalem]{ulem}

\newtheorem{theorem}{Theorem}[section]
\newtheorem{lemma}[theorem]{Lemma}

\theoremstyle{definition}
\newtheorem{definition}[theorem]{Definition}
\newtheorem{example}[theorem]{Example}

\newtheorem{prop}[theorem]{Proposition}

\theoremstyle{plain}
\newtheorem{cor}[theorem]{Corollary}
\newtheorem{proposition}[theorem]{Proposition}

\allowdisplaybreaks
\theoremstyle{remark}
\newtheorem{remark}[theorem]{Remark}

\allowdisplaybreaks

\newcommand{\E}{\calE}

\newcommand{\dist}{\mathrm{dist}}
\newcommand{\diam}{\operatorname{diam}}

\newcommand{\mres}{\mathbin{\vrule height 1.6ex depth 0pt width
0.13ex\vrule height 0.13ex depth 0pt width 1.3ex}}
\newcommand*{\dd}{\mathop{}\!\mathrm{d}}

\newtheorem{step}{Step}

\makeatletter
\@addtoreset{step}{proposition}
\@addtoreset{step}{theorem}
\makeatother

\makeatletter
\newcommand*\bigcdot{\mathpalette\bigcdot@{.5}}
\newcommand*\bigcdot@[2]{\mathbin{\vcenter{\hbox{\scalebox{#2}{$\m@th#1\bullet$}}}}}
\makeatother

\newcommand{\R}{\mathbb{R}}
\renewcommand{\S}{\mathbb{S}}
\newcommand{\sphere}{\mathbb{S}}

\newcommand{\N}{\mathbb{N}}

\newcommand{\Div}{\operatorname{div}}

\newcommand{\calE}{\mathcal{E}}

\newcommand{\calH}{\mathcal{H}}
\newcommand{\calM}{\mathcal{M}}
\renewcommand{\H}{\mathcal{H}}
\newcommand{\wto}{\rightharpoonup}

\newcommand{\spt}{\mathrm{spt}}
\newcommand{\IVo}{\mathsf{IV}^o}
\renewcommand{\dd}{\,\mathrm{d}}
\DeclareMathOperator{\sgn}{sgn}

\def\opt{\mathrm{w}}

\makeatletter
\@namedef{subjclassname@2020}{%
  \textup{2020} Mathematics Subject Classification}
\makeatother

\numberwithin{equation}{section}

\begin{document}
\title[Elastic bilayer membranes under confinement]{Elastic bilayer membranes under confinement \\ --- Existence, regularity, and rigidity}

\author[M.~Röger]{Matthias Röger}
\address[M.~Röger]{Biomathematics group, Department of Mathematics, TU Dortmund University, Vogelpothsweg 87, 44227 Dortmund, Germany.}
\email{matthias.roeger@tu-dortmund.de}

\author[F.~Rupp]{Fabian Rupp}
\address[F.~Rupp]{Faculty of Mathematics, University of Vienna, Oskar-Morgenstern-Platz 1, 1090 Vienna, Austria.}
\email{fabian.rupp@univie.ac.at}




\date{\today}

\begin{abstract}
Motivated by applications to cell biology, we study the constrained minimization of the Helfrich energy among closed surfaces confined to a container. 
We show existence of minimizers in the class of bubble trees of spherical weak branched immersions and derive the Euler--Lagrange equations which involve a measure-valued Lagrange multiplier that is concentrated on the coincidence set with the container boundary. We provide a careful analysis of this elliptic system and prove optimal regularity for solutions throughout the branch points. For surfaces confined in the unit ball we show that the minimization problem behaves rather rigid and identify a parameter range for which minimizers are always round spheres.
\end{abstract}
\maketitle


\section{Introduction}
Cells are the fundamental building blocks of living organisms.
Their functions are closely related to their shapes and can vary widely, depending on a multitude of factors, both internal and external.
Most cells are separated from their surroundings by elastic membranes, typically in form of \emph{lipid bilayers}.
Since these membranes are extremely thin, they are commonly modeled by two-dimensional structures, i.e., surfaces.
Among the most prominent theoretical frameworks is the \emph{Canham--Helfrich model} \cites{Canham,Helfrich} which characterizes the equilibrium shapes of such membranes as (local) minimizers of a curvature-dependent bending energy under various constraints.
In particular, this model provides a mechanical explanation of the characteristic biconcave shape of red blood cells \cite{DeulingHelfrich}.

For a smooth oriented surface $S\subset \R^3$, the energy considered by Helfrich \cite{Helfrich} takes the form
\begin{align}
    \frac{\kappa_H}{4}\int_{S}(H-H_0)^2\dd A + \kappa_K \int_S K\dd A\, .
\end{align} Here $K$ is the Gauss curvature, $H$ is the trace of the scalar second fundamental form of $S$ with respect to the unit normal field $\nu$ along $S$ determined by the orientation, and $\dd A$ denotes integration with respect to the area measure. The constants $\kappa_H, \kappa_K,H_0 \in\R$ are modeling parameters. If the topology is fixed, the Gauss curvature integral is already determined as a consequence of the Gauss--Bonnet Theorem. We thus consider $\kappa_K=0$ and $\kappa_H=1$ in the sequel and  the  \emph{(reduced) Helfrich energy}
\begin{align}\label{eq:intro_def_Helfrich}
    \calE(S) \vcentcolon = \calE_{H_0}(S)\vcentcolon = \frac{1}{4}\int_{S}(H-H_0)^2\dd A\,.
\end{align}
Our convention is that when $S$ is embedded, we assume that $\nu$ is the interior unit normal, yielding $H = \vec{H}\cdot \nu \equiv \frac{2}{r}$ for $S=\partial B_r(0)$ a sphere of radius $r>0$.
The remaining parameter $H_0\in\R$ in \eqref{eq:intro_def_Helfrich} is called the \emph{spontaneous curvature} and accounts for a preferred curvature of the membrane.
To incorporate the inextensibility of the membrane or the impermeability, variational models often involve constraints of prescribed area $|S|=a$ for some $a>0$, or prescribed volume.

In the particular case of zero spontaneous curvature $H_0=0$, we obtain the \emph{Willmore energy} $\mathcal{W}\vcentcolon=\calE_0$.
This functional has many additional distinct properties, most importantly its conformal invariance \cite{blaschke}.
By Willmore's inequality $\mathcal{W}\geq 4\pi$ for all closed immersed surfaces, with equality only for embedded round spheres \cite{Willmore65}.

In this work, we study the Helfrich energy \eqref{eq:intro_def_Helfrich} for closed surfaces in the presence of an \emph{obstacle constraint} as a reduced model for vesicles that are confined by a \emph{container}.
This problem is motivated by biological applications, for example elastic organelles trapped by a rather rigid cell wall (for instance, mitochondriae inside a plant cell), or red blood cells moving through capillaries, which can be smaller in diameter than the cells themselves.

\subsection{Previous work}
The rigorous mathematical treatment of the variational Canham--Helfrich model started with the work of Schygulla \cite{MR2928137} who showed the existence of smooth minimizers with spherical topology, assuming $H_0=0$ and a constraint on both area and enclosed volume.
This has recently been extended to the higher genus case in \cites{KellerMondinoRiviere2014,Scharrer22NLA,MondinoScharrer2023,KusnerMcGrath23}.
For $H_0\neq 0$, Mondino--Scharrer \cite{MR4076069} showed existence of spherical constrained minimizers in the class of \emph{bubble trees of weak branched immersions} \cite{MR3276119}.
This framework is generically necessary, since bubbling phenomena may occur and yield nonsmooth minimizers \cite{RuppScharrer23}*{Example 1.2}.
Moreover, existence results in weak Geometric Measure Theory settings (without prescribing the topology) have been obtained in \cites{BrazdaLussardiStefanelli20,KubinLussardiMorandotti24}.

Formally, the (scalar) Euler--Lagrange equation of the constrained Helfrich energy takes the form
\begin{align}\label{eq:intro_ELE}
        \Delta_g H + (\frac12 H^2- 2K) H + 2 H_0 K -(2\Lambda+\frac12 H_0^2) H - 2\rho = 0\,,
\end{align}
where $\Delta_g$ is the Laplace--Beltrami operator and $\Lambda,\rho\in \R$ are some Lagrange multipliers arising from the area and volume constraint, respectively.
Already for $H_0=0$, \eqref{eq:intro_ELE} involves cubic terms in the curvature that do not need to be integrable for  immersions in Sobolev spaces having only finite Willmore energy.
Remarkably, using local conformal coordinates away from the branch points, Rivi\`ere \cite{MR2430975} provided a formulation of \eqref{eq:intro_ELE} for weak immersions with finite energy and, even more, showed that it is equivalent to a system of conservation laws with a \emph{Jacobian structure.} This was extended in \cite{MR4076069} to the case $H_0\neq 0$ and essentially used in the discussion of regularity of minimizers.
Moreover, in \cite{MR4706029}, these techniques have recently been adapted to also allow for a general class of inhomogeneities in \eqref{eq:intro_ELE}.
The behavior and regularity at potential branch points is more complicated and, in the case of the Willmore energy, has been investigated in several works \cites{KuSc04,MR2318282,MR3096502,MR4076072}.
The weak formulation of the Euler--Lagrange equations, away from finitely many branch points, and its Jacobian structure is a key advantage of weak immersions (and, more generally, bubble trees) over other weak notions of surfaces such as (curvature) varifolds, where such concepts are not known, except in a small energy regime \cite{rupp2024global}.

Confinement constraints in variational models for Helfrich-type energies have been proposed and analyzed in several contributions from biophysics. 
In \cites{KaSM12,KaSM12a}, spherical inclusions were considered and Euler--Lagrange equations in rotational symmetry were derived alongside numerical simulations.
In further works, additional adhesion energies were taken into account \cite{RiPK14}, and minimization in classes of toroidal vesicles was explored \cite{BoMM15}.

A mathematical analysis of the minimization of the Willmore energy under a spherical confinement condition and for fixed surface area was presented in \cite{MuRo14}. 
A lower bound estimate for the Willmore energy of confined topological spheres in terms of the prescribed area was proven, along with a sharp increase in energy when the area slightly exceeds that of the container. 
The lower bound was extended \cite{Wojt17} to classes of surfaces with general topology, leading to various consequences for minimization problems involving Helfrich energies.

The variational problem for differences $\mathcal{W}-\Lambda\mathcal{A}$ between Willmore and area functionals under a general confinement condition has been analyzed in \cite{Pozz23}, characterizing the values of the parameter $\Lambda$, for which the minimal energy is finite and proving existence and regularity properties of minimizers.

Additionally, a minimization problem for the Willmore functional concerning surfaces of revolution under Dirichlet (clamped) boundary conditions and an obstacle constraint has been addressed in \cite{GrOk23}. 
For obstacles where the energy infimum falls below a suitable threshold value, existence and regularity of minimizers have been established.

Finally, it is worth mentioning that numerous contributions have considered the minimization of the elastica functional for curves---a one-dimensional analogue to the Willmore energy---such as those found in \cites{DoMR11,DaDe18,DaMN18,Muel19,Wojt21,GrOk23} and their references.

\subsection{Main contributions}

\subsubsection{Existence}
Our first contribution is an existence result for the confinement problem for the Helfrich energy in an arbitrary container without any a priori assumption on the energy infimum.
The most obvious class to search for minimizers might be the set of embedded or immersed spheres.
However, due to the aforementioned bubbling phenomena \cite{RuppScharrer23}*{Example 1.2}, this is not a particularly promising variational framework, see also \Cref{ex:sphere_packing} and \Cref{rem:sphere_packing} below.
Instead, as in \cite{MondinoScharrer2023}, we rely on the theory of bubble trees and their compactness properties developed in \cite{MR3276119}. Intuitively, a bubble tree can be viewed as the union of finitely many spherical $W^{2,2}$-immersions such that the parametrizations are glued together in suitable sense, see \Cref{subsec:bubble_trees} for the relevant definitions.

\begin{theorem}[See \Cref{thm:existence_bubble} below]\label{thm:intro:existence}
    Let $\Omega\subset\R^3$ be nonempty open set.
    Let $a>0$ and $H_0\in\R$.
    Then the Helfrich energy $\E_{H_0}$ attains its minimum in the class of bubble trees of weak branched immersions with prescribed area $a$ that are contained in $\overline{\Omega}$.
\end{theorem}

\subsubsection{Regularity} Our second focus in this work is on the regularity properties of minimizers.
Under a mild nondegeneracy assumption on the area constraint, we construct suitable test functions that allow us to derive the weak Euler--Lagrange equations in a local conformal parametrization, see \Cref{lem:measure-ELE} below.
We show that this equation is satisfied \emph{throughout the branch points} and involves a \emph{measure-valued inhomogeneity} due to the confinement constraint, in contrast to previous works.
In the spirit of \cite{MR2430975}, we show that the Euler--Lagrange equations can be transformed into a system of conservation laws with a Jacobian structure.
This is not immediate, since we allow for branch points, and thus the previous derivations in  \cites{MR2430975,MR3524220}, see also \cites{MR3518329,MR4076069}, 
do not directly apply.
\begin{theorem}[See \Cref{thm:reg_global} below]\label{thm:intro:regularity_global}
    Let $\Omega$ be of class $C^3$.
    Then every minimizer in \Cref{thm:intro:existence} consists of bubbles that are $W^{2,p}$-regular for all $p\in [1,\infty)$, in particular $C^{1,\alpha}$ for all $\alpha<1$, even at the branch points.
    Moreover, away from the container boundary and the branch points, the bubbles are smooth.
\end{theorem}

In fact, we show $W^{3,p}$-regularity of solutions to the weak Euler--Lagrange equations for all $p\in[1,2)$, see \Cref{thm:reg_meas_ELE-2}. By considering the \emph{inverted catenoid,} this regularity result is optimal, in general, see \Cref{lem:inv_cat} below. 

Since our local regularity analysis allows for a measure term, we can also treat \emph{branch point singularities} (which generate Dirac contributions), even in the absence of confinement conditions. 
In \Cref{sec:futher_applications}, we gather some implications of our work for the regularity through the branch points of (constrained) minimizers, including those studied by Mondino--Scharrer \cite{MondinoScharrer2023} and Da Lio--Palmurella--Rivi\`ere \cite{MR4076072}, see \Cref{cor:MS_regularity} and \Cref{cor:DLPR_regularity} below. Moreover, our analysis applies to \emph{branched Willmore immersions.} It has been shown in \cites{KuSc04,MR2318282} that the optimal (Hölder) regularity in this case is $C^{1,\alpha}$ for all $\alpha<1$. Our findings provide the optimal regularity in the third order Sobolev scale, see \Cref{cor:branch_Willmore_regularity} below. 
\subsubsection{Rigidity} 
The scaling properties of the geometric quantities imply that the energy of a round sphere $\partial B_r(x_0)$ is given by
\begin{align}\label{eq:energy_sphere}
    \calE_{H_0}(\partial B_r(x_0))= (2/r-H_0)^2 \pi r^2= (2-H_0 r)^2\pi\,.
\end{align}
As in \cite{MuRo14}, we now consider the particular case where the container $\Omega=B\vcentcolon = B_1(0)$ is the unit ball.
First, we do not impose an area constraint and minimize $\E$ among smooth embedded closed surfaces contained in $\overline{B}$.
If $H_0$ favors small spheres, it is easy to see that the minimal energy vanishes.
More surprisingly, the preference for spheres is rather rigid and extends to the regime $0\leq H_0\leq 2$. The following scaling argument shows that outside of the range $0\leq H_0\leq 2$, the unit sphere can not be optimal. 
For any $S$ and $t>0$, we have
\begin{equation}
    4\calE_{H_0}(tS) -4\calE_{H_0}(S)
    = 2(1-t)H_0\int_{S} (H-H_0)\dd A
    +(1-t)^2H_0^2|S|\,.
\label{eq:scaling_arg}
\end{equation}
Therefore, $S$ can only be minimal if
\begin{equation}
    H_0\int_{S} (H-H_0)\dd A\geq 0\,.
\end{equation}
In particular, the unit sphere is a candidate for minimizing $\E_{H_0}$ only if $0\leq H_0\leq 2$. Our third main contribution shows that it is indeed the minimizer in this parameter range.
\begin{theorem}[See \Cref{thm:main_ext_obstacle} below]\label{thm:intro:main_ext_obstacle}
Let $0\leq H_0\leq 2$ and let $S$ be a smooth, embedded closed surface with $S\subset \overline{B}$.
Then we have
    \begin{align}\label{eq:intro:E(B)_minimal}
        \E_{H_0}(S)\geq \calE_{H_0}(\partial B)
        =(2-H_0)^2\pi\,.
    \end{align}
    Equality is attained if and only if $S=\partial B$.
\end{theorem}
The rigidity result depends crucially on the value of $H_0\in\R$ (see \Cref{thm:main_ext_obstacle}) and can be generalized in various ways.
We first extend it to \emph{varifolds} (\Cref{prop:rigidity_varifolds}) and then to the area constrained confinement problem as considered in \Cref{thm:intro:existence} (\Cref{prop:rigidity_area_constraint}).
Moreover, we characterize the dependence of the infimum bending energy on the fixed area $a>0$ and the value of the spontaneous curvature $H_0\in\R$, see \Cref{pro:est-impr-below}.
As in \cite{MuRo14}*{Theorem 3}, we also show that, in the regime $H_0\in [0,2)$, the infimum grows sharply like a square root as the prescribed area $a$ just exceeds $4\pi$, see \Cref{rem:sharp_square_root}.

\section{Preliminaries}
In this section, we recall the notation and concepts that we use for our variational study and provide, for the convenience of the reader, an overview of the relevant results.

\subsection{Notation}\label{sec:notation}
To avoid using indices and to shorten notation, we adapt the conventions from \cite{MondinoScharrer2023} based on \cite{MR3524220}. We denote by $\cdot$ the inner product in $\R^n$ and use $\times$ for the cross product in $\R^3$. Let $D=B_1(0)\subset \R^2$ be the unit disk. 
For (in a suitable sense) differentiable maps $\Phi, \Psi\colon D\to\R^3$, we write
\begin{align}
    \nabla \Phi &\vcentcolon = \begin{pmatrix}
        \partial_{x^1} \Phi \\
        \partial_{x^2}\Phi\\
    \end{pmatrix}\, , & \nabla^\perp \Phi &\vcentcolon=
    \begin{pmatrix}
        -\partial_{x^2} \Phi \\
        \partial_{x^1}\Phi\\
    \end{pmatrix}\, ,\\
     \langle \Phi, \nabla \Psi\rangle &\vcentcolon = \begin{pmatrix}
        \Phi\cdot \partial_{x^1} \Psi \\
        \Phi \cdot \partial_{x^2}\Psi\\
    \end{pmatrix}\, , &
    \langle \Phi, \nabla^\perp \Psi\rangle &\vcentcolon 
    = \begin{pmatrix}
        - \Phi \cdot \partial_{x^2}\Psi \\
        \Phi \cdot \partial_{x^1}\Psi\\
    \end{pmatrix}\, ,\\ 
    \Phi\times \nabla \Psi &\vcentcolon = \begin{pmatrix}
        \Phi\times \partial_{x^1}\Psi \\
        \Phi\times \partial_{x^2}\Psi\\
    \end{pmatrix}\, , & 
    \Phi\times \nabla^\perp \Psi &\vcentcolon = \begin{pmatrix}
        -\Phi\times \partial_{x^2}\Psi \\
        \Phi\times \partial_{x^1}\Psi\\
    \end{pmatrix}\,,
\end{align}
and
\begin{align}
      \nabla \Phi\times \nabla \Psi &\vcentcolon = 
        \partial_{x^1}\Phi\times \partial_{x^1}\Psi + \partial_{x^2}\Phi\times \partial_{x^2}\Psi\, , & 
        \nabla \Phi\times \nabla^\perp \Psi&\vcentcolon = -\partial_{x^1}\Phi\times \partial_{x^2}\Psi + \partial_{x^2}\Phi\times \partial_{x^1}\Psi\, , \\
        \nabla \Phi \cdot \nabla \Psi &\vcentcolon = \partial_{x^1}\Phi\cdot\partial_{x^1}\Psi + \partial_{x^2}\Phi\cdot \partial_{x^2}\Psi\,, &   \nabla \Phi \cdot \nabla^\perp \Psi &\vcentcolon = -\partial_{x^1}\Phi\cdot\partial_{x^2}\Psi + \partial_{x^2}\Phi\cdot \partial_{x^1}\Psi\, .
\end{align}
Note that with this convention, $\cdot$ always yields a scalar, whereas $\langle\cdot,\cdot\rangle$ and $\times$ yield vectors or tuples of vectors. Moreover, for $\lambda\colon D\to\R$, define
\begin{align}
    \langle \nabla \lambda, \Phi\rangle &\vcentcolon = \begin{pmatrix}
        \partial_{x^1}\lambda\,\Phi \\
        \partial_{x^2}\lambda\,\Phi\\
    \end{pmatrix}\, , &
    \langle \nabla^\perp \lambda, \Phi\rangle&\vcentcolon = \begin{pmatrix}
        -\partial_{x^2}\lambda\, \Phi \\
        \partial_{x^1}\lambda\, \Phi
    \end{pmatrix}\, ,\\
    \langle \nabla \lambda, \nabla\Phi\rangle &\vcentcolon = \langle \nabla^\perp\lambda, \nabla^\perp\Phi\rangle & \langle \nabla^\perp \lambda, \nabla\Phi\rangle &\vcentcolon = -\langle \nabla\lambda, \nabla^\perp\Phi\rangle  \, \\
    &\vcentcolon=\partial_{x^1}\lambda\partial_{x^1}\Phi + \partial_{x^2}\lambda \partial_{x^2}\Phi\,, &    
    &\vcentcolon = 
        -\partial_{x^2}\lambda \partial_{x^1}\Phi +\partial_{x^1}\lambda \partial_{x^2}\Phi\, .
\end{align}
Lastly, for $X = \begin{pmatrix}
    X_1 \\X_2
\end{pmatrix}$, where $X_1,X_2\colon D\to\R^3$, the divergence is
\begin{align}
    \Div X \vcentcolon = \partial_{x^1}X_1 + \partial_{x^2}X_2\,,
\end{align}
so that in particular the Laplacian is given by $\Delta \Phi = \Div (\nabla \Phi)$.

\subsection{Oriented varifolds}\label{sec:OrVarif}
As one concept of generalized surfaces that allow for a notion of curvature and a generalization of the Helfrich functional, we use oriented varifolds as introduced by Hutchinson \cite{Hutc86}, see also \cite{Allard}.
For convenience we restrict ourselves to the case of $2$-varifolds in $\R^3$.
We identify the set of $2$-dimensional oriented subspaces of $\R^3$ with the unit vectors $\sphere^2\subset\R^3$.
With this identification, an \emph{oriented $2$-varifold in $\R^3$} is a Radon measure on $\R^3\times \sphere^2$.

\begin{definition}
    Consider a $2$-rectifiable set $S\subset\R^3$, an orientation on $S$ given by an $\calH^2$-measurable unit normal field $\nu:S\to\sphere^2$, and locally $\calH^2$-summable multiplicity functions $\theta_\pm:S\to\N_0$ with $\theta_++\theta_->0$.
    We associate to the tuple $(S,\nu,\theta_+,\theta_-)$ the oriented $2$-varifold
    \begin{align}
        \underline{v}(S,\nu,\theta_+, \theta_-)(\varphi)
        =\int_{S\times \sphere^2}\big(\theta_+(x)\varphi(x,\nu(x))+\theta_-(x)\varphi(x,-\nu(x))\big)\dd \calH^2(x)\,.
    \end{align}
    The set of oriented $2$-varifolds $V$ that can be represented in the form $\underline{v}(S,\nu,\theta_+,\theta_-)$ defines the class of \emph{integral oriented $2$-varifolds} $\IVo_2(\R^3)$.
    The Radon measure $\mu=\Vert V\Vert = (\theta_++\theta_-)\calH^2\mres S$ over $\R^3$ is called the \emph{mass of $V$}.
    We say that $V$ has \emph{weak mean curvature} $\vec H$ if there exists $\vec H\in L^1_{\mathrm{loc}}(\mu;\R^3)$ such that
    \begin{align}
        \int \Div_\mu X(x) \dd \mu(x) = -\int \vec H(x)\cdot X(x)\dd \mu(x)\qquad \text{ for all }X\in C_c^1(\R^3;\R^3)\,,
    \end{align}
    where, for $\mu$-a.e.\ $x\in \R^3$, $\Div_\mu X(x)$ denotes the tangential divergence of $X$ in $x$ with respect to the approximate tangent plane $T_x\mu$ of $\mu$ at $x$.
\end{definition}

Any smooth proper immersion of an oriented surface without boundary naturally induces a smooth unit normal field $\nu$ and a varifold $V\in \IVo_2(\R^3)$ with weak mean curvature $\vec H\in L^2(\mu;\R^3)$, see \cite{RuppScharrer23}*{Section 2.3}.
In this case, $\vec H = H \nu$ is the classical mean curvature with $H=\kappa_1+\kappa_2$.
Even for a general $V=\underline{v}(S,\nu,\theta_+,\theta_-)\in \IVo_2(\R^3)$ we have that $\vec H(x)$ and $\nu(x)$ are parallel for $\mu$-a.e.\ $x\in \R^3$, see \cite{Brakke}*{Section 5.8}.

\begin{definition}[Helfrich energy for oriented varifolds]
\label{def:HelfOrVar}
Consider any integral oriented $2$-varifold $V=\underline{v}(S,\nu,\theta_+,\theta_-)$ with weak mean curvature $\vec H\in L^2(\mu;\R^3)$, where $\mu=\|V\|$.
For $H_0\in\R$ we define the generalized Helfrich energy of $V$ by
\begin{equation}
    \calE_{H_0}(V) :=
    \int_{S} |\vec{H}-H_0\nu|^2\theta_+ + |\vec{H}+H_0\nu|^2\theta_- \dd\calH^2\,. \label{eq:HelfOrVar}
\end{equation}
\end{definition}

This definition is consistent with \cite{MR3896647}*{Section 2} and \cite{RuppScharrer23}*{Section 2.2}.
Note that reversing the orientation of the varifold corresponds to switching the sign of $H_0$ in the Helfrich energy, similarly for immersions, cf.\ \cite{RuppScharrer23}*{(1.1)}.
If we want to remove this degree of freedom and give an interpretation to the sign of $H_0$, we need to adapt the convention  for the unit normal in the embedded case described just after equation \eqref{eq:intro_def_Helfrich} above.
This naturally leads to varifolds satisfying a suitable \emph{divergence theorem}: the \emph{varifolds with enclosed volume} \cite{RuppScharrer23}*{Section 4.3} and \emph{volume varifolds} \cite{scharrer2023properties}*{Definition 3.1}.
In the smooth case, immersions with an appropriate notion of ``interior'' are \emph{Alexandrov immersion} \cite{MR143162}, see \cite{RuppScharrer23}*{Definition 1.4}.
They correspond to the boundary values of an immersion of a compact $3$-manifold with boundary into $\R^3$.
The key property is that, for all these classes, a certain \emph{concentrated volume} introduced in \cite{RuppScharrer23} is positive at any given point, see \cite{RuppScharrer23}*{Lemma 4.9} and \cite{RuppScharrer23}*{(5.1)}.

\subsection{Bubble trees of weak branched immersions}\label{subsec:bubble_trees}

In order to derive existence results for the minimization of the Canham-Helfrich energy under a confinement constraint, we will use the notion of weak branched immersion, which was formalized in \cites{MR3059839,MR2928715,MR3276154} and which we recall in the following.

\begin{definition}\label{def:weak_imm}
Let $\Sigma$ be an oriented surface and let $g_0$ be a smooth reference metric on $\Sigma$.
We call a map $\Phi\in W^{2,2}\cap W^{1,\infty}(\Sigma;\R^3)$ a \emph{weak branched immersion} if and only if
\begin{enumerate}
    \item there exists $1<C<\infty$ such that $C^{-1} |\mathrm{d}\Phi|_{g_0}^2 \leq |\mathrm{d}\Phi \wedge \dd \Phi|_{g_0}\leq C|\mathrm{d}\Phi|_{g_0}^2$, where in local coordinates $\dd \Phi \wedge\dd \Phi\vcentcolon = (\dd x^1\wedge \dd x^2) \partial_{x^1}\Phi\wedge \partial_{x^2} \Phi$;
    \item there is a finite (possibly empty) set of \emph{branch points} $b_1,\dots, b_N\in \Sigma$, $N\in \N_0$, such that
    \begin{align}
        \log |\mathrm{d}\Phi|_{g_0}\in L^\infty_{\mathrm{loc}}(\Sigma\setminus \{b_1,\dots, b_N\})\,;
    \end{align}
    \item the \emph{Gauss map} $\nu$, defined in a positive chart $x$ by $$\nu\vcentcolon = \frac{\partial_{x^1} \Phi \times \partial_{x^2} \Phi}{|\partial_{x^1} \Phi \times \partial_{x^2} \Phi|}\,,$$
    satisfies $\nu\in W^{1,2}(\Sigma;\R^3)$.
\end{enumerate}
If $N=0$ (i.e., in the absence of branch points), we call $\Phi$ a \emph{($W^{2,2}$-)Lipschitz immersion.} The space of weak branched immersions does not depend on the choice of $g_0$ and is denoted by $\mathcal{F}_\Sigma$.

Moreover, we say that $\Phi$ is \emph{conformal} (with respect to $g_0$) if in all conformal charts $x$ of $(\Sigma,g_0)$, we have $|\partial_{x^1}\Phi| = |\partial_{x^2}\Phi|$ and $\partial_{x^1}\Phi \cdot \partial_{x^2}\Phi=0$.
In this case, we call $\lambda\vcentcolon = \log|\mathrm{d}\Phi|_{g_0}$ the \emph{conformal factor.}  
We will mostly work with the sphere $\Sigma=\S^2$ with the standard round metric or with the disk $\Sigma=D=\{z\in \R^2: |z|<1\}$ with the flat metric.
\end{definition}

Conformal immersions play a key role in the analysis of the Willmore energy. It is well-known consequence of the work of Müller--\v{S}ver\'ak \cite{MR1366547}, see \cite{MR2928715}, that a conformal parametrization may always be arranged locally for immersions in $\mathcal{F}_\Sigma$, see also \cite{MR3524220} for an approach based on H\'elein's moving frame technique \cite{MR1913803}.

\begin{lemma}\label{lem:exist_conf_repara}
	Let $\Phi\in \mathcal{F}_\Sigma$ and $x\in \Sigma$. Then there exists an open neighborhood $U\subset\Sigma$ of $x$ and a bilipschitz diffeomorphism $\varphi\colon D\to U$ such that $\Phi\circ \varphi \in \mathcal{F}_D$ is conformal.
\end{lemma} 

Fixing local positive coordinates $x^1,x^2$ on $\Sigma$ and denoting by $g=(g_{ij})_{ij}$ the pullback metric of the Euclidean inner product under $\Phi \in \mathcal{F}_\Sigma$, the \emph{mean curvature} and the \emph{mean curvature vector} are defined by
\begin{align}
    H \vcentcolon= g^{ij} \partial_{x^i}\partial_{x^j} \Phi \cdot\nu = - g^{ij} \partial_{x^j} \Phi\cdot\partial_{x^i} \nu\,, \qquad \vec H\vcentcolon = H \nu\,,
\end{align}
respectively.
In local coordinates, the area measure induced by the metric $g$ is given by $\mu = \sqrt{\det g} \dd x^1\otimes \dd x^2$.
We may define the area and Helfrich energy of $\Phi\in \mathcal{F}_\Sigma$ by
\begin{align}
     \mathcal{A}(\Phi) \vcentcolon = \int_\Sigma \dd \mu\,,\quad \E(\Phi) \vcentcolon = \frac{1}{4}\int_\Sigma (H-H_0)^2 \dd \mu\,.
\end{align}

\Cref{def:weak_imm} leads to the definition of \emph{bubble trees}, see \cite{MR3276119}*{Definition 7.1}.
\begin{definition}\label{def:bubble_tree}
    We call an $(N+1)$-tuple $ T = (f,\Phi^1,\dots,\Phi^N)$, $N\in\N$, a \emph{bubble tree of (spherical) weak branched immersions} if $f\in W^{1,\infty}(\S^2;\R^3)$ and $\Phi^1, \dots,\Phi^N\in \mathcal{F}_{\S^2}$ are  such that the following holds.
    \begin{itemize}
        \item There exist open geodesic balls $B^1, \dots, B^N\subset \S^2$ such that $\overline{B^1}=\S^2$ and for all $i\neq i'$ either $\overline{B^i}\subset B^{i'}$ or $\overline{B^{i'}}\subset B^i$.
        \item For all $i\in\{1,\dots,N\}$ there exists $N^i\in \N$ and disjoint geodesic balls $B^{i,1},\dots,B^{i,N^i}\subset\S^2$ with $\overline{B^{i,1}},\dots,\overline{B^{i,N^i}}\subset B^i$ such that for all $i'\neq i$ either $\overline{B^i}\subset B^{i'}$ or $\overline{B^{i'}}\subset B^{i,j}$ for some $j\in\{1,\dots,N^i\}$.
        \item For all $i \in \{1, \ldots, N\}$ there exist distinct points $b^{i,1}, \ldots, b^{i,N^i} \in \mathbb S^2$ and a Lipschitz diffeomorphism
	\begin{align}
		\Xi^i: B^i \setminus \bigcup_{j = 1}^{N^i - 1} \overline{B^{i,j}} \to \mathbb S^2 \setminus \{b^{i,1}, \ldots, b^{i, N^i}\}
	\end{align}
	with Lipschitz extension
	\begin{align}
		\overline \Xi_i: \overline{B^i} \setminus \bigcup_{j=1}^{N^i - 1}B^{i,j} \to \mathbb S^2
	\end{align}
	such that
	\begin{align}
		\overline{\Xi}_i(\partial B^{i,j}) = b^{i,j} \text{ whenever } j \in \{1, \ldots, N^{i}-1\}\,, \qquad \overline \Xi_i(\partial B^i) = b^{i,N^i}.
	\end{align}
    \item For all $i \in \{1, \ldots, N\}$,
	\begin{align} \label{bt:parametrisation}
		f(x) = (\Phi^i \circ \Xi^i)(x) \qquad \text{whenever } x \in B^i \setminus \bigcup_{j = 1}^{N^{i}-1} \overline{B^{i,j}}
	\end{align}
	and for all $j \in \{1, \ldots, N^i\}$ there exists $p^{i,j} \in \mathbb R^3$ such that
	\begin{align} \label{bt:branch_point}
		f(x) = p^{i,j} \qquad \text{whenever } x \in B^{i,j} \setminus \bigcup_{i' \in J^{i,j}} \overline{B^{i'}}
	\end{align}
	where $J^{i,j} = \{i': \overline{B^{i'}} \subset B^{i,j}\}$.
    \end{itemize}
    We refer to $\Phi^1,\dots,\Phi^N$ as \emph{bubbles} of $T$.
\end{definition}
For brevity, we will occasionally refer to a bubble tree of weak branched immersions simply as a \emph{bubble tree}. Note that if $T=(f, \Phi^1, \dots,\Phi^N)$ is  a bubble tree, then
    \begin{align}\label{eq:image_f_phi^i}
        f(\S^2) = \bigcup_{i=1}^N \Phi^i(\S^2)\,,
    \end{align}
so $T$ really corresponds to a union of weak branched immersions. We thus define
    \begin{align}
        \diam  T \vcentcolon = \diam f(\S^2)\,,\qquad
        \calE_{H_0} (T)\vcentcolon = \sum_{i=1}^N \calE_{H_0}(\Phi^i)\,,\qquad \mathcal{A}( T)\vcentcolon = \sum_{i=1}^N \mathcal{A}(\Phi^i)\,.
    \end{align}

The space of bubble trees satisfies the following essential compactness result.

\begin{theorem}[{\cite{MR4076069}*{Theorem 3.3}, see also \cite{MR3276119}*{Theorem 7.2}}]\label{thm:compactness}
    Let $ T_k = (f_k, \Phi_k^1,\dots, \Phi_k^{N_k})$ be a sequence of bubble trees of weak branched immersions with
    \begin{align}
        \limsup_{k\to\infty} \sum_{i=1}^{N_k} \int_{\S^2} (1+|\vec H_{\Phi_k^i}|^2)\dd\mu_{\Phi_k^i} <\infty\,, \qquad \liminf_{k\to\infty}\sum_{i=1}^{N_k}\diam \Phi_k^i(\S^2)>0\,.
    \end{align}
    Then, after passing to a subsequence, we have $N_k\equiv N\in \N$ for all $k\in\N$ and there exist diffeomorphisms $\Psi_k$ of $\S^2$ with
    \begin{align}\label{eq:bubble_C^0-conv}
        f_k\circ \Psi_k\to  u_\infty\qquad \text{ in }C^0(\S^2;\R^3)\,,
    \end{align}
    for some $ u_\infty\in W^{1,\infty}(\S^2;\R^3)$.
    Moreover, for all $i\in\{1,\dots,N\}$ there exists $Q^i\in\N$ and sequences $\Psi_k^{i,1},\dots,\Psi_k^{i,Q^i}$ of positive conformal diffeomorphisms of $\S^2$ such that for each $j\in \{1,\dots, Q^i\}$ there exist finitely many points $b^{i,j,1},\dots,b^{i,j,Q^{i,j}}\in\S^2$ with
    \begin{align}
        \Phi_k^i\circ \Psi_k^{i,j}\wto \xi^{i,j}_\infty \qquad \text{ in }W^{2,2}_{\mathrm{loc}}(\S^2\setminus \{b^{i,j,1},\dots,b^{i,j,Q^{i,j}}\};\R^3)
    \end{align}
    for some conformal weak branched immersion $\xi^{i,j}_\infty\in \mathcal{F}_{\S^2}$.
    Furthermore,
    \begin{align}
         T_\infty\vcentcolon = ( u_\infty, (\xi_\infty^{1,j})_{j=1,\dots, Q^1},\dots, (\xi^{N,j}_\infty)_{j=1,\dots,Q^N})
    \end{align}
    is a bubble tree of weak branched immersions and we have
    \begin{align}
        \mathcal{A}( T_\infty)=\lim_{k\to\infty}\mathcal{A}( T_k)\,, \qquad \calE_{H_0}(T_\infty)\leq \liminf_{k\to\infty}\calE_{H_0}(T_k)\,.
    \end{align}
\end{theorem}

\subsection{Localized first variation}\label{sec:prelims_variation}

For $\Phi\in \mathcal{F}_\Sigma$ and $\zeta\in L^\infty(\Sigma)$, we define a localized version of the Helfrich energy by
\begin{align}
    \calE(\Phi;\zeta) \vcentcolon = \frac{1}{4}\int_\Sigma (H-H_0)^2\zeta\dd \mu\,,
\end{align}  
and analogously, we define $\mathcal{A}(\Phi; \zeta)$.
Further, if $\Phi\in \mathcal{F}_D$ conformal and $\zeta\in L^\infty(\Sigma)$, we define the distributions
    \begin{align}
        \delta\calE(\Phi;\zeta).w\vcentcolon
        =& \int_D \frac{1}{2} (\vec H-H_0\nu)\cdot \Delta w\,\zeta \dd z  \\
        &\quad + \int_D\frac{1}{4}\left((3 H - 4H_0)\nabla \nu - \vec H \times \nabla^\perp \nu -(2H_0 H -H_0^2)\nabla \Phi\right)\cdot\nabla w\,\zeta\dd z\,,\label{eq:def_dE}\\
        \delta\mathcal{A}(\Phi;\zeta).w \vcentcolon =& \int_D \nabla \Phi\cdot \nabla w \,\zeta \dd z\,,\label{eq:def_dA}
    \end{align}
acting on $w\in C^\infty(D;\R^3)$. If $\zeta\equiv 1$, we will just write $\delta\calE(\Phi)$ and $\delta\mathcal{A}(\Phi)$, respectively. The notation is motivated by the fact that, for $\Phi$ conformal, $\delta\E(\Phi)$ coincides with the first variation of $\E$ for test functions supported away from the branch points, see \cite{MR3524220}*{Theorem 4.57}, \cite{MR4076069}*{Lemma 4.1}, and also \Cref{lem:local_first_vari}. 

Now, let $\Sigma$ be a closed oriented surface and let $\Phi\in \mathcal{F}_\Sigma$. We extend $\delta\calE$ and $\delta\mathcal{A}$ to $\Phi\in \mathcal{F}_{\Sigma}$ and $w\in C^\infty(\Sigma;\R^3)$. First, around any point $x\in \Sigma$, by \Cref{lem:exist_conf_repara} there exists an open set $U$ and a bilipschitz homeomorphisms $\varphi\colon D\to U$ such that $\Phi\circ \varphi\in \mathcal{F}_D$ is conformal. For $\zeta\in L^\infty(U)$ and $w\in C^\infty(U;\R^3)$, we define
\begin{align}\label{eq:def_dE_chart}
    \delta\calE(\Phi;\zeta).w \vcentcolon = \delta\calE(\Phi\circ\varphi;\zeta\circ \varphi).(w\circ\varphi)\,,
\end{align}
using formula \eqref{eq:def_dE}. This is well-defined, i.e., independent of the choice of chart $\varphi$ above. Indeed, the transition map between any such charts is a conformal diffeomorphisms of the disk. A direct computation using the transformation of the integrand in \eqref{eq:def_dE} under composition with conformal diffeomorphisms from the right yields the claim.

Since $\Sigma$ is compact, it can be covered by the images of finitely many charts $\varphi_i\colon D\to U_i$, $i\in I$, as above. Let $\zeta_i \in C_c^\infty(U_i)$, $i\in I$, be a partition of unity subordinate to this cover and define
\begin{align}\label{eq:def_dE_partition_of_unity}
    \delta\calE(\Phi).w \vcentcolon = \sum_i \delta\calE(\Phi;\zeta_i).w \quad \text{for }w\in C^\infty(\Sigma;\R^3)\, .
\end{align}
Using linearity in $\zeta$ of \eqref{eq:def_dE} and \eqref{eq:def_dE_chart}, it follows that this does not depend on the particular partition of unity used. 

For $\mathcal{A}$ one may proceed exactly in the same way. The notation $\delta \calE(\Phi), \delta\mathcal{A}(\Phi)$ is justified by the following observation which we will prove in \Cref{sec:local_variaton}.

\begin{lemma}\label{lem:E_A_diff}
    Let $\Sigma$ be a closed oriented surface and let $\Phi\in \mathcal{F}_\Sigma$ with branch points $\{b_1,\dots,b_N\}$.
    The functionals $\calE$ and $\mathcal{A}$ are Gateaux differentiable in $\Phi\in \mathcal{F}_{\Sigma}$ with respect to perturbations  $w\in W^{2,2}\cap W^{1,\infty}(\Sigma;\R^3)$ with compact support in $\Sigma\setminus\{b_1,\dots,b_N\}$ and
    \begin{align}\label{eq:E_A_diff_1}
        \delta\calE(\Phi).w =  \left.\frac{\dd}{\dd t}\right\vert_{t=0}\calE(\Phi+tw)\,,\qquad
        \delta\mathcal{A}(\Phi).w = \left.\frac{\dd}{\dd t}\right\vert_{t=0}\mathcal{A}(\Phi+tw)\,.
    \end{align}
\end{lemma}

Following the arguments in \cite{MR2989995}*{Lemma A.5}, one can show that the functionals $\calE$, $\mathcal{A}$ are in fact Fr\'echet differentiable. However, for our purposes the Gateaux differentiability in \eqref{eq:E_A_diff_1} suffices.

\section{Minimizers for the general confinement problem}
We now study confinement problems with prescribed area in a general container.
A natural space to obtain compactness is the one of bubble trees of weak branched immersions, see \Cref{subsec:bubble_trees}.
In this section we consider an arbitrarily fixed value $H_0\in\R$ and for simplicity write $\calE$ instead of $\calE_{H_0}$.

\subsection{Existence in the class of bubble trees}

We have the following existence result.

\begin{theorem}\label{thm:existence_bubble}
    Let $\Omega\subset \R^3$ be a nonempty open set.
    Let $a>0$, let $H_0\in\R$, and let
    \begin{align}
        \mathcal{M}(a,\overline{\Omega})\vcentcolon=\{  T  = (f, \Phi^1, \dots, \Phi^N)\text{ bubble tree of weak branched immersions} \mid  \mathcal{A}( T) = a, f(\S^2)\subset \overline{\Omega}\}\,.
    \end{align}
    Then, there exists $ T_*\in\mathcal{M}(a,\overline{\Omega})$ such that
    \begin{align}\label{eq:min_prob_bubble_tree}
        \calE( T_*) = \min_{ T\in \mathcal{M}(a,\overline{\Omega})}\calE( T)\,.
    \end{align}
\end{theorem}

\begin{proof}
    First, note that $\mathcal{M}(a,\overline{\Omega})\neq\emptyset$.
    Indeed, taking sufficiently many copies of a small sphere connected with catenoidal necks (see \Cref{ex:sphere_packing} below) will satisfy the assumptions.
    Let
    \begin{align}
    \inf_{ T\in \mathcal{M}(a,\overline{\Omega})}\calE( T)=\vcentcolon E<\infty\,.
    \end{align}
    Let $ T_k=(f_k,\Phi^1_k,\dots,\Phi^{N_k}_k)\in \mathcal{M}(a,\overline{\Omega})$ be a minimizing sequence.
    Then the Willmore energy of each $\Phi_k^i$ is bounded since, for $k$ sufficiently large,
    \begin{align}\label{eq:existence_willmore_bound}
        \mathcal{W}(\Phi_k^i)\leq \mathcal{W}( T_k) \leq 2 \calE( T_k)+ 2 H_0^2 a \leq 2 (E+1)+ 2H_0^2 a\,.
    \end{align}
    Applying Simon's lower diameter estimate \cite{Simon}*{Lemma 1.1} to each $\Phi^i_k$, we thus obtain
    \begin{align}
        \diam \Phi^i_k(\S^2)\geq \sqrt{\frac{\mathcal{A}(\Phi^i_k)}{2(E+1)+ 2H_0^2 a}}\,.
    \end{align}
    Using the inequality $\sqrt{x+y}\leq \sqrt{x}+\sqrt{y}$ for $x,y\geq 0$, we find
    \begin{align}
        \sum_{i=1}^{N^i_k}\diam\Phi^i_k(\S^2)\geq \sqrt{\frac{\sum_{i=1}^{N^i_k}\mathcal{A}(\Phi^i_k)}{2 (E+1) + 2H_0^2 a}} =  \sqrt{\frac{{a}}{2 (E+1) + 2H_0^2 a}}\,.
    \end{align}
    We may thus apply \Cref{thm:compactness} and, after passing to a subsequence, obtain a limit bubble tree $ T_\infty$.
    With the convergence \eqref{eq:bubble_C^0-conv}, we find that for each $x\in\S^2$ we have
    \begin{align}
         u_\infty (x) = \lim_{k\to\infty} f_k(\Psi_k(x)) \in \overline{\Omega}\,,
    \end{align}
    since $f_k(\S^2)\subset \overline{\Omega}$.
    Moreover, it follows from \Cref{thm:compactness} that $\mathcal{A}( T_\infty)=a$, so $ T_\infty\in\mathcal{M}(a,\overline{\Omega})$.
    Together with the lower semicontinuity of $\calE$ in \Cref{thm:compactness}, we conclude
    \begin{align}
        \calE( T_*) = \inf_{ T\in\mathcal{M}(a,\overline{\Omega})}\calE( T)\,. &\qedhere
    \end{align}
\end{proof}

\begin{remark}\label{rem:existence_volume}
    As in \cite{MR4076069}, one may also impose an additional volume constraint in the variational problem \eqref{eq:min_prob_bubble_tree}.
    For any $\Phi\in\mathcal{F}_{\Sigma}$, the (algebraic) volume is given by
    \begin{align}
        \mathcal{V}(\Phi)\vcentcolon = -\frac{1}{3}\int_\Sigma \nu\cdot \Phi \dd \mu\,,
    \end{align}
    and for a bubble tree as in \Cref{def:bubble_tree}, we may define $\mathcal{V}( T) \vcentcolon = \sum_{i=1}^N \mathcal{V}(\Phi^i)$.
    Since the constraint $\mathcal{V}( T)=v$ is closed under the convergence of bubble trees by \cite{MR4076069}*{Lemma 3.1}, one also has existence of minimizers for $\calE$ in the class of bubble trees $ T$ with $\mathcal{A}( T)=a, \mathcal{V}( T)=v, f(\S^2)\subset\overline{\Omega}$, provided the admissible class is not empty. Similarly, as it was done in \cite{MR4865096}, one could impose an (additional) constraint on the total mean curvature
    \begin{align}
        \mathcal{T}(\Phi)\vcentcolon = \int_\Sigma H\dd \mu\,,
    \end{align}
    which is continuous under the convergence of bubble trees \cite{MondinoScharrer2023}*{Theorem 3.3}.    
    Since neither of this is  necessary to have a well-posed variational problem, we will only consider an area constraint in the sequel.
    However, for possible future applications, in \Cref{subsec:conservation_laws} below we also allow for Lagrange multipliers corresponding to the volume and the total mean curvature.
\end{remark}

\begin{remark}
Similarly as in \cite{RuppScharrer23}*{Section 6.3}, one may estimate the concentrated volume defined in \cite{RuppScharrer23} for \emph{embedded} surfaces of prescribed area confined in $\overline{\Omega}$ by their Helfrich energy. Thus, assuming the infimum energy to lie below a certain threshold, the Hausdorff density along a minimizing sequence will be strictly below $2$. Following the arguments in \cite{RuppScharrer23}*{Section 6.3}, this may be used to  solve the minimization problem \eqref{eq:min_prob_bubble_tree} in the class of \emph{Lipschitz quasi-embeddings} \cite{RuppScharrer23}*{p.~34}. Identifying this threshold explicitly and constructing appropriate comparison surfaces is an interesting task for future research.
\end{remark}

We conclude this section with a diameter bound on minimizers of \eqref{eq:min_prob_bubble_tree}.

\begin{lemma}\label{lem:diam_bound}
    Any minimizer in \Cref{thm:existence_bubble} has diameter bounded by
    \begin{align}
        \diam  T_* \leq \frac{1}{\pi} \Big(2 \min_{ T\in \mathcal{M}(a,\overline{\Omega})}\calE( T) + (2H_0^2+1)a\Big)\,.
    \end{align}
\end{lemma}
\begin{proof}
    Note that $f(\S^2)$ is connected by continuity.
    Hence, if $ T_*= (f, \Phi^1,\dots,\Phi^N)$, by \eqref{eq:image_f_phi^i} we have
    \begin{align}
        \diam  T_* \leq \sum_{i=1}^N \diam \Phi^i(\S^2)\leq \frac{2}{\pi} \sum_{i=1}^N \sqrt{\mathcal{W}(\Phi^i) \mathcal{A}(\Phi^i)} \leq \frac{1}{\pi}\Big(\mathcal{W}( T)+\mathcal{A}( T)\Big)\,,
    \end{align}
    where, for each $\Phi^i$, we used Simon's upper diameter inequality \cite{Simon} with the constant $C=2/\pi$, see \cite{MR1650335}*{Lemma 1}.
    The statement follows from the area constraint and \eqref{eq:existence_willmore_bound}.
\end{proof}

\subsection{Derivation of the Euler--Lagrange equations}

 In this section, we will derive the Euler--Lagrange equations for the variational problem \eqref{eq:min_prob_bubble_tree}.
 We now assume that $\Omega$ is $C^3$-regular.
  Let $n$ be the interior unit normal along $\partial\Omega$, and, for $\tau>0$, define
\begin{align}
    \Omega_\tau = \{ y+ tn(y)\mid y\in \partial\Omega, |t|<\tau \}\,.
\end{align}
Note that by \Cref{lem:diam_bound}, any minimizer in \Cref{thm:existence_bubble} has a priori bounded diameter.
 Hence, even for $\Omega$ unbounded, we may assume that $\partial\Omega$ possesses a tubular neighborhood of uniform size $\tau>0$ in the sequel in which the projection $\Pi\colon \Omega_\tau\to \partial\Omega, x+tn(x)\mapsto x$ is of class $C^2$ and any $y\in \Omega_\tau$ can be uniquely written as
\begin{align}\label{eq:tub_nbd}
    y= \Pi(y) + d(y) n(y)\,,
\end{align}
where $d\in C^3(\Omega_\tau)$ is the signed distance function and $n$ is extended to a $C^2$-function on all of $\Omega_\tau$ via $n(y) = n(\Pi(y))$.
Note that with a representation as
in \eqref{eq:tub_nbd}, we have $y\in \overline{\Omega} \cap \Omega_\tau$ if and only if $d(y) \geq 0$.
 
As a starting point for our regularity analysis, we fix a minimizer $ T_* = (f,\Phi^1,\dots,\Phi^N)$ in \Cref{thm:existence_bubble} and discuss the regularity of some fixed $\Phi^i\in \mathcal{F}_{\S^2}$ with branch points $b^{i,1},\dots,b^{i,N^i}$.
 Since \eqref{eq:min_prob_bubble_tree} is of constrained obstacle type, we will need a nondegeneracy assumption on the area constraint.
 Hence in the following we assume that
 \begin{align}\label{eq:constr_nondeg}
     \text{ there exists }i_0\in\{1,\dots,N\} \text{ with }\int_{(\Phi^{i_0})^{-1}(\Omega)}|\vec H_{\Phi^{i_0}}|\dd\mu_{\Phi^{i_0}} >0\,.
 \end{align}
 
 Our goal in this section is to find the weak Euler--Lagrange equations for the bubbles $\Phi^i$.
 \begin{lemma}\label{lem:measure-ELE}
    Assume \eqref{eq:constr_nondeg} and let $i\in \{1,\dots,N\}$.
    There exists a vectorial Radon measure $\vec \sigma$ over $\S^2\setminus\{b^{i,1},\dots,b^{i,N^i}\}$, $ a^{i,1},\dots,  a^{i,N^i}\in\R^3$, and $\Lambda\in\R$ such that we have
	\begin{align}\label{eq:measure-ELE}
        \delta\calE(\Phi^i) + \Lambda\delta\mathcal{A}(\Phi^i) =  \vec \sigma + \sum_{k=1}^{N^i} a^{i,k} \delta_{b^{i,k}} \quad \text{ in }\mathcal{D}'(\S^2;\R^3)\,.
	\end{align}
    Moreover, $\spt\vec \sigma\subset (\Phi^i)^{-1}(\partial\Omega)\setminus\{b^{i,1},\dots,b^{i,N^i}\}$ and $\vec \sigma = n\circ\Phi^i \sigma$ for some nonnegative Radon measure $\sigma$ over $\S^2\setminus\{b^{i,1},\dots,b^{i,N^i}\}$.
\end{lemma}

\begin{proof}
    For fixed $i\in\{1,\dots,N\}$, we abbreviate $\Phi\vcentcolon = \Phi^i$ and write $B \vcentcolon= \{b^{i,1},\dots,b^{i,N^i}\}$ for the branch points.

    \begin{step}\label{step:ele-1}
        We consider $w\in W^{2,2}\cap W^{1,\infty}(U,\R^3)$ with $U\subset\subset \S^2\setminus B$ such that $\Phi(\overline{U})\subset \Omega$.
     Suppose first that \eqref{eq:constr_nondeg} is satisfied with $i_0=i$.
     It follows from \eqref{eq:constr_nondeg} and the continuity of $\Phi$ that there exists an  open set $V\subset \S^2\setminus B$ with $\Phi(\overline{V})\subset \Omega$, $\int_{V}|\vec H_{\Phi}|\dd\mu_{\Phi}>0$.
     Hence, there exists $ u\in C_c^\infty(\S^2;\R^3)$ with $\spt u\subset V$ and, by the first variation of area and \Cref{lem:E_A_diff},
    \begin{align}
        \delta\mathcal{A}(\Phi). u = -\int \vec H_{\Phi}\cdot  u \dd\mu_{\Phi} >0\,.
    \end{align}
   We may then define the variation
    \begin{align}
        \Phi_{t,s} \vcentcolon = \Phi + t w  + s  u\,.
    \end{align}
    For all $t\in (-\varepsilon,\varepsilon)$, $s\in (-\varepsilon,\varepsilon)$ with $\varepsilon>0$ sufficiently small, using that $\dd \Phi$ only vanishes in $B$, we have that $\Phi_{t,s}$ is well-defined, $\Phi(0,0)=\Phi$, and $\Phi_{t,s}\in \mathcal{F}_{\S^2}$.
    By \Cref{lem:E_A_diff}, the map $(t,s)\mapsto \mathcal{A}(\Phi_{t,s})$ is $C^1$ and satisfies
    \begin{align}\label{eq:lem:constr_nondeg_i=i0}
        \left.\frac{\partial}{\partial s} \mathcal{A}(\Phi_{t,s})\right\vert_{(t,s)=(0,0)} = \delta\mathcal{A}(\Phi). u >0\,.
    \end{align}
    Reducing $\varepsilon>0$, the Implicit Function Theorem yields the existence of a $C^1$-map $\alpha\colon(-\varepsilon,\varepsilon)\to \R$ such that $\alpha(0)=0$ and such that with
    \begin{align}
        \Phi_t \vcentcolon = \Phi_{t,\alpha(t)}\,,
    \end{align}
    we have $\mathcal{A}(\Phi_t) = \mathcal{A}(\Phi)$ for all $t\in (-\varepsilon,\varepsilon)$.
    Since $\dist(\Phi(\overline{U}), \partial\Omega)>0$ and $\dist(\Phi(\overline{V}),\partial\Omega)>0$, we conclude that with $|t|<\varepsilon$ with $\varepsilon>0$ sufficiently small, we have $\Phi_{t}(\S^2)\subset \overline{\Omega}$.
    Consequently, for any  $t\in (-\varepsilon,\varepsilon)$, we replace $\Phi$ in $T_*$ by $\Phi_t$ and, after suitably modifying $f$, still have an admissible bubble tree in $\mathcal{M}(a,\overline{\Omega})$.
    By \Cref{lem:E_A_diff}, $\calE$ is differentiable  with respect to $W^{2,2}\cap W^{1,\infty}$-perturbations away from the branch points, so by minimality of $ T_*$ we find that
    \begin{align}
        0= \left.\frac{\dd}{\dd t}\right\vert_{t=0}\calE(\Phi_t) = \delta\calE(\Phi).w + \alpha'(0)\delta\calE(\Phi).  u\,.
    \end{align}
   On the other hand, by \Cref{lem:E_A_diff} and the choice of $\alpha$ we have
    \begin{align}\label{eq:dtA=0}
        0 = \left.\frac{\dd}{\dd t}\right\vert_{t=0}\mathcal{A}(\Phi_t) = \delta \mathcal{A}(\Phi). w + \alpha'(0) \delta \mathcal{A}(\Phi). u \,.
    \end{align}
    Choosing
    \begin{align}\label{eq:def_lambda_i=i0}
    \Lambda \vcentcolon= - \frac{\delta\calE(\Phi). u}{ \delta\mathcal{A}(\Phi). u} \in\R\,,
    \end{align}
    cf.\ \eqref{eq:lem:constr_nondeg_i=i0}, we find
    \begin{align}\label{eq:vari_ineq_step_1}
        0 = \delta\calE(\Phi).w + \Lambda\delta\mathcal{A}(\Phi).w\,.
    \end{align}
    
    If \eqref{eq:constr_nondeg} is satisfied for some $i_0\neq i$, writing $\Psi \vcentcolon =\Phi^{i_0}$, we find an open set $W\subset \S^2\setminus\{b^{i_0,1},\dots,b^{i_0,N^{i_0}}\}$ with $\Psi(\overline{W})\subset \Omega$, and a function
    $v\in C_c^\infty(W;\R^3)$ such that
    \begin{align}\label{eq:lem:constr_nondeg}
        \delta\mathcal{A}(\Psi).v >0\,.
    \end{align}
    Similar as above, for $\varepsilon>0$ small there exists a $C^1$-function $\beta\colon(-\varepsilon,\varepsilon)\to\R$ with $\beta(0)=0$ such that, for $t\in (-\varepsilon,\varepsilon)$
    \begin{align}
        \Phi_t \vcentcolon = \Phi + tw\,,\quad
        \Psi_t \vcentcolon = \Psi + \beta(t)v\label{eq:def_Psi_t}
    \end{align}
    are well-defined, $\Phi_t \in\mathcal{F}_{\S^2},\Psi_t\in\mathcal{F}_{\S^2}$, $\Phi_t(\S^2)\cup\Psi_t(\S^2)\subset\overline{\Omega}$, and $\mathcal{A}(\Phi_t)+\mathcal{A}(\Psi_t)=\mathcal{A}(\Phi)+\mathcal{A}(\Psi)$.
    In this case, in $ T_*$ we may replace $\Phi$ by $\Phi_t$and $\Psi$ by $\Psi_t$, yielding again (after appropriately adjusting $f$) an admissible bubble tree in $\mathcal{M}(a,\overline{\Omega})$.
    Using \Cref{lem:E_A_diff}, the minimality of $ T_*$, the area constraint, and \eqref{eq:constr_nondeg} we conclude as above that with
    \begin{align}\label{eq:def_lambda_other}
        \Lambda \vcentcolon = - \frac{\delta\calE(\Psi).v}{\delta\mathcal{A}(\Psi).v}\in\R\,,
    \end{align}
   we still have \eqref{eq:vari_ineq_step_1}. $\hfill\diamond$
    \end{step}
    
    \begin{step}\label{step:ele-2}
        We prove the statement for $w = \xi n\circ\Phi$, where $\xi\in W^{2,2}\cap W^{1,\infty}(U)$ with $U\subset\subset \S^2\setminus B$ and $\Phi(\overline U)\subset \Omega_\tau$.
        Suppose first that $\xi\geq 0$ and that \eqref{eq:constr_nondeg} is satisfied for $i=i_0$.
        Let $ u$, $V$ be as in \Cref{step:ele-1} and consider the variation
    \begin{align}
        \Phi_{t,s} \vcentcolon = \Phi + t w + s  u = \Pi\circ \Phi + (d\circ\Phi + t\xi)n\circ\Phi+s u\,.
    \end{align}
    As above, for $|t|,|s|<\varepsilon$ small, we have that $\Phi_{t,s}$ is well-defined, $\Phi(0,0)=\Phi$, $\Phi_{t,s}\in \mathcal{F}_{\S^2}$, since $n,d,\Pi$ are of class $C^2$.
    Again, due to \eqref{eq:lem:constr_nondeg_i=i0}, there exists a $C^1$-map $\alpha\colon(-\varepsilon,\varepsilon)\to \R$ such that $\alpha(0)=0$ and
    \begin{align}
        \Phi_t \vcentcolon = \Phi_{t,\alpha(t)}
    \end{align}
    satisfies $\mathcal{A}(\Phi_t) = \mathcal{A}(\Phi)$ for all $t\in (-\varepsilon,\varepsilon)$.
    We consider $t\in [0,\varepsilon)$.
    In $V$, by the choice of $ u$ in \Cref{step:ele-1}, we have $\dist(\Phi(\overline{V}),\partial\Omega)>0$, so $\Phi_{t}(\overline{V})\subset \Omega$ for $|t|<\varepsilon$ sufficiently small.
    In $\S^2\setminus V$, $ u\equiv 0$, so that by \eqref{eq:tub_nbd} and since $\xi\geq 0$, we have $\Phi_{t}(\S^2)\subset \overline{\Omega}$.
    Replacing $\Phi$ by $\Phi_t$ as in \Cref{step:ele-1} --- now only for $t\in[0,\varepsilon)$ --- and using minimality of $ T_*$ yields
    \begin{align}
        0\leq \left.\frac{\dd}{\dd t}\right\vert_{t=0}\calE(\Phi_t) = \delta\calE(\Phi).w + \alpha'(0)\delta\calE(\Phi).  u\, .
    \end{align}
    On the other hand, as above
    \begin{align}
        0 = \left.\frac{\dd}{\dd t}\right\vert_{t=0}\mathcal{A}(\Phi_t) = \delta \mathcal{A}(\Phi). w + \alpha'(0) \delta \mathcal{A}(\Phi). u\, .
    \end{align}
    Choosing $\Lambda$ as in \eqref{eq:def_lambda_i=i0} yields
    \begin{align}\label{eq:vari_ineq_xi}
        0\leq \delta\calE(\Phi).w + \Lambda\delta\mathcal{A}(\Phi).w\, .
    \end{align}
      
    If \eqref{eq:constr_nondeg} is satisfied for some $i_0\neq i$, we may argue similarly, constructing $\Phi_t, \Psi_t$ using the same $v$ as in \Cref{step:ele-1}, to conclude that \eqref{eq:vari_ineq_xi} holds
    for $\Lambda$ as in \eqref{eq:def_lambda_other}.
        
    Now \eqref{eq:vari_ineq_xi} implies that $L\in \mathcal{D}'(\Phi^{-1}(\Omega_\tau)\setminus B)$ defined by
    \begin{align}
        \langle L, \xi\rangle \vcentcolon = \delta\calE(\Phi).(\xi n\circ\Phi) + \Lambda \delta\mathcal{A}(\Phi).(\xi n\circ\Phi)\,, \qquad \xi\in C_c^\infty(\Phi^{-1}(\Omega_\tau)\setminus B)\,,
    \end{align}
    is a nonnegative distribution.
    Hence, by the Riesz representation theorem, there exists a nonnegative Radon measure $\sigma$ over $\Phi^{-1}(\Omega_\tau)\setminus B$ representing $L$.
    As a consequence, $\langle L, \xi\rangle = \int \xi\dd\sigma$ for all $\xi\in W^{2,2}\cap W^{1,\infty}(U)$ with $U\subset\subset \S^2\setminus B$ and $\Phi(\overline U)\subset\Omega_\tau$.
    Using \Cref{step:ele-1}, i.e., Equation \eqref{eq:vari_ineq_step_1}, for $\xi$ compactly supported in $\Phi^{-1}(\Omega_\tau\cap \Omega)$ we conclude that $\sigma$ is supported in $\Phi^{-1}(\partial\Omega)\setminus B$.
    $\hfill\diamond$
    \end{step}

    \begin{step}\label{step:ele-3}
    We prove the statement for $w\in W^{2,2}\cap W^{1,\infty}(U;\R^3)$ with $U\subset\subset\S^2\setminus B$, $\Phi(\overline{U})\subset \Omega_{\tau/2}$, $\Vert w\Vert_\infty^2 \leq \tau/2$, and $ w\cdot n\circ \Phi =0$.
    
    Suppose that \eqref{eq:constr_nondeg} holds with $i_0=i$.
    Take $U'$ open such that $\overline U\subset U'\subset \S^2\setminus B$ and $\Phi(\overline{U'})\subset \Omega_{\tau/2}$ and let $\chi \in C_c^\infty(\S^2)$with $\chi_{U} \leq \chi \leq \chi_{U'}$.
    With $ u$ and $V$ as in \Cref{step:ele-1}, we define
    \begin{align}
        \Phi_{t,s} \vcentcolon = \Pi\circ\Phi + (d\circ\Phi + t^2\chi )n\circ\Phi + tw+s u= \Phi + t^2\chi n\circ\Phi + tw + s  u\, .
    \end{align}
    We have $\Phi_{0,0}=\Phi$ and, for all $t,s\in (-\varepsilon,\varepsilon)$ with $\varepsilon>0$ sufficiently small, we have that $\Phi_{t,s}\in \mathcal{F}_{\S^2}$.
    Arguing as in \Cref{step:ele-1}, \eqref{eq:lem:constr_nondeg_i=i0} and the Implicit Function Theorem yield the existence of $\alpha\colon(-\varepsilon,\varepsilon)\to\R$ with $\alpha(0)=0$ so that for all $t\in (-\varepsilon,\varepsilon)$ we have $\Phi_t\vcentcolon = \Phi_{t,\alpha(t)}\in \mathcal{F}_{\S^2}$ and $\mathcal{A}(\Phi_t)=\mathcal{A}(\Phi)$.  Since $\dist(\Phi(\overline V),\partial\Omega)>0$ for $|t|<\varepsilon$ with $\varepsilon>0$ sufficiently small, we have $\Phi_{t}(\overline{V})\subset \overline{\Omega}$.
    By \eqref{eq:tub_nbd}, we have that
    \begin{align}
        y_t \vcentcolon = \Pi\circ\Phi + (\tau/2+ d\circ\Phi)n\circ\Phi
    \end{align}
    satisfies $y_t(x)\in \Omega\cap \Omega_\tau$  for $x\in \overline U$
    with $\dist(y_t(x), \partial\Omega)=\tau/2+d\circ\Phi(x)$.
    On the other hand, in $\overline U\setminus \overline V$ we have $ u\equiv 0$, $\chi\equiv 1$, so, recalling $\Vert w\Vert_\infty^2 \leq \tau/2$, we have
    \begin{align}
        |\Phi_{t} - y_t|^2 = (t^2-\tau/2)^2 + t^2|w|^2 \leq (\tau/2)^2 \quad \text{ for }|t|^2\leq \tau/2\,.
    \end{align}
    Thus, by the definition of $y_t$, we conclude that $\Phi_{t}(\overline U)\subset \overline{\Omega}$ provided $|t|<\varepsilon$ for some $\varepsilon>0$ sufficiently small.
    Lastly, in $U'\setminus \overline{U} \subset \Phi^{-1}(\Omega_{\tau/2})$, we have $\Phi_{t} = \Phi + t^2 \chi n\circ\Phi \in \overline{\Omega}\cap \Omega_\tau$ by \eqref{eq:tub_nbd} as $|t|^2\leq \tau/2$.
    Overall, we thus have that $\Phi_{t}(\S^2)\subset \overline{\Omega}$ for $|t|\leq \varepsilon$ small enough.
    
    As in \Cref{step:ele-1}, minimality of $T_*$ and the choice of $\alpha$ yield
    \begin{align}
        0&= \left.\frac{\dd}{\dd t}\right\vert_{t=0}\calE(\Phi_t) = \delta\calE(\Phi).w + \alpha'(0)\delta\calE(\Phi).  u\,,\\
        0 &= \left.\frac{\dd}{\dd t}\right\vert_{t=0}\mathcal{A}(\Phi_t) = \delta \mathcal{A}(\Phi). w + \alpha'(0) \delta \mathcal{A}(\Phi). u\,,
    \end{align}
    so that with exactly the same $\Lambda$ as in \eqref{eq:def_lambda_i=i0}, we have
    \begin{align}
         0= \delta\calE(\Phi).w + \Lambda\delta\mathcal{A}(\Phi).w\,.
    \end{align}

    The case where \eqref{eq:constr_nondeg} is satisfied for some $i_0\neq i$ can be treated similarly, constructing $\Phi_t$, $\Psi_t$ using the same $v$ as in \Cref{step:ele-1}, hence yielding the same $\Lambda$ as in \eqref{eq:def_lambda_other}.$\hfill\diamond$
    \end{step}

    \begin{step}\label{step:ele-4}
        Now consider $w \in W^{2,2}\cap W^{1,\infty}(U;\R^3)$ with $U\subset\subset \S^2\setminus B$.
        First assume, in addition, that $\Vert w \Vert_\infty^2 \leq \tau/2$.
        Let $\eta$ be a smooth cutoff function with $\eta\equiv 1$ in $\Phi^{-1}(\Omega_{\tau/4})$ and $\eta \equiv 0$ in $\Phi^{-1}(\Omega\setminus \Omega_{\tau/2})$.
        Set $w_1 \vcentcolon= (1-\eta) w$, $w_2\vcentcolon =  \eta w$.
        Further, we may decompose
    $w_2 = \xi n\circ\Phi + r$
    with $\xi \vcentcolon= w_2 \cdot n\circ \Phi  \in W^{2,2}\cap W^{1,\infty}(U)$ and $r\in W^{2,2}\cap W^{1,\infty}(U;\R^3)$ such that $r\cdot  n\circ\Phi  =0$.
    Since the same $\Lambda\in \R$ is used in Steps \ref{step:ele-1}--\ref{step:ele-3}, cf.\ \eqref{eq:def_lambda_i=i0} and \eqref{eq:def_lambda_other}, we have
    \begin{align}
        \delta\calE(\Phi).w_1 + \Lambda \delta\mathcal{A}(\Phi).w_1 &=0\,,\\
        \delta\calE(\Phi).(\xi n \circ \Phi) + \Lambda \delta\mathcal{A}(\Phi).(\xi n \circ \Phi) &= \langle \xi, \sigma\rangle\,, \\
        \delta\calE(\Phi).r + \Lambda \delta\mathcal{A}(\Phi).r &= 0\,.
    \end{align}
    By linearity, \eqref{eq:measure-ELE} is proven for all $w\in C_c^\infty(\S^2\setminus B;\R^3)$ with $\vec \sigma\vcentcolon = n\circ\Phi\sigma \mres \Phi^{-1}(\Omega_\tau)$. $\hfill\diamond$
    \end{step}

    \begin{step}\label{step:ele-5}
        Consequently, the distribution $ L\in \mathcal{D}'(\S^2;\R^3)$ given by
     \begin{align}
        \langle  L,w\rangle \vcentcolon =\delta\calE(\Phi).w + \Lambda \delta\mathcal{A}(\Phi).w -\int_{\S^2\setminus B}  w\cdot n\circ \Phi \dd \sigma\,, \quad w\in C^\infty(\S^2;\R^3)\,,\label{eq:def_L}
    \end{align}
    is supported in the finite set $B$.
    It follows, using for instance \cite{MR1996773}*{Theorem 2.3.4}, that $ L$ is a weighted sum of derivatives of the Dirac masses $\delta_b$ with $b\in B$.
    Choose a local conformal chart $\xi\colon D\to\S^2$ in a small neighborhood of $b\in B$ with $\xi(0)=b$.
    We may use the definition of $\delta\calE, \delta\mathcal{A}$, see \eqref{eq:def_dE} and \eqref{eq:def_dA}, to conclude that for all $w\in C_c^\infty(D;\R^3)$ we have
    \begin{align}
        \langle  L,w\rangle &= \frac{1}{2}\int_D  (\vec H-H_0\nu) \cdot \Delta w\dd z + \int_D a\cdot \nabla w \dd z - \int_{D\setminus \{0\}} w\cdot n\circ\Phi \dd \sigma \\
        &= \sum_{|\alpha|\leq k}  a_\alpha\cdot \nabla^\alpha w(0)\,.\label{eq:represent_L}
    \end{align}
    Here the sum runs over multi-indices $\alpha$ with norm $|\alpha|\leq k<\infty$, $ a_\alpha \in \R^3$ are coefficients, and $ a \in L^1(D;(\R^3)^2)$ by \eqref{eq:def_dE} and \eqref{eq:def_dA} since $\Phi\in \mathcal{F}_D$.
    We now test \eqref{eq:represent_L} with $ w_\varepsilon(z):= w(z/\varepsilon)$, $z\in D$, $\varepsilon>0$, and observe that
    \begin{align}
        &\frac{\varepsilon^{-2}}{2} \int_{D} \vec H \cdot \Delta w(\cdot/\varepsilon) \dd z + \varepsilon^{-1} \int a\cdot \nabla w(\cdot/\varepsilon)\rangle \dd z + \int_{D\setminus \{0\}}  w(\cdot/\varepsilon) \cdot n\circ\Phi \dd \sigma \\
        &= \sum_\alpha \varepsilon^{-|\alpha|}  a_\alpha\cdot \nabla^\alpha w(0)\,.\label{eq:no_delta_derivative_1}
    \end{align}
    By Hölder's inequality and change of variables, we may estimate the modulus of the left hand side by
    \begin{align}
        \frac{\varepsilon^{-1}}{2} \Vert \vec H \Vert_{L^2(B_\varepsilon(0))} \Vert \Delta w\Vert_{L^2(D)} + \varepsilon^{-1} \Vert  a \Vert_{L^1(B_\varepsilon(0))} \Vert \nabla w\Vert_{L^\infty(D)} + \Vert w\Vert_\infty \sigma(D\setminus \{0\})\,.\label{eq:no_delta_derivative_2}
    \end{align}
    Since the vector $(\nabla^\alpha w(0))_{|\alpha|\leq k})$ can take any arbitrary value for a suitable choice of $w\in C_c^\infty(D)$, together with the behavior of \eqref{eq:no_delta_derivative_1} and \eqref{eq:no_delta_derivative_2} as $\varepsilon\to 0$, one deduces that $ a_\alpha = 0$ for all $|\alpha|\geq 1$.
    Proceeding like this around every branch point $b\in B$, we deduce $ L = \sum_{k=1}^{N^i}  a^{i,k} \delta_{b^{i,k}}$ for some $ a^{i,k}\in\R^3$ as in the statement.\qedhere
    \end{step}
\end{proof}

We now show that the Dirac contribution vanishes at branch points that do not touch the boundary $\partial\Omega$.
The key ingredient is the first variation formula in \Cref{lem:local_first_vari} for variations with noncompact support and a suitable vector field, adapted to cancel the Dirac weight, see \cite{KuSc04}*{p.~338} and also \cite{MR4076072}*{Lemma 4.11}.

\begin{lemma}\label{lem:residue}
    In the setting of \Cref{lem:measure-ELE}, for any $i=1,\dots,N$, $k=1,\dots,N^i$, if $\Phi(b^{i,k})\in\Omega$, we have $ a^{i,k}=0$.
\end{lemma}

\begin{proof}
    Fix $i\in \{1,\dots,N\}$ and write $\Phi\vcentcolon = \Phi^i$, $B \vcentcolon = \{b^{i,1},\dots,b^{i,N^i}\}$.
    Let $b\vcentcolon=b^{i,k}\in B$ and let $ a\vcentcolon =  a^{i,k}$ be the weight of the corresponding Dirac measure in \Cref{lem:measure-ELE}.
    
    Suppose first that \eqref{eq:constr_nondeg} is satisfied with $i_0=i$.
    Let $ u\in C_c^\infty(\S^2\setminus B;\R^3)$ be as in \Cref{step:ele-1} of the proof of \Cref{lem:measure-ELE}.
    By continuity, there exist open sets $b\in U_1\subset \subset U_2$ such that $\overline{U_2}\cap B=\{b\}$, $\Phi(\overline{U_2})\subset \Omega$, and $\overline{U_2}\cap \spt  u =\emptyset$.
    Thus, there exists $w\in C_c^\infty(\S^2;\R^3)$ such that
    \begin{align}
        w = \begin{cases}
             a & \text{in } U_1\,,\\
            0 & \text{in } \S^2\setminus \overline{U_2}\,.
        \end{cases}
    \end{align}
    Proceeding as in \Cref{step:ele-1} of the proof of \Cref{lem:measure-ELE}, we find that there exists $\varepsilon>0$, $\alpha\in C^1((-\varepsilon,\varepsilon);\R)$ with $\alpha(0)=0$ such that $\Phi_t \vcentcolon = \Phi+ tw +\alpha(t) u \in \mathcal{F}_{\S^2}$, $\Phi_t (\S^2)\subset\overline{\Omega}$, and $\mathcal{A}(\Phi_t)=\mathcal{A}(\Phi)$ for $t\in (-\varepsilon,\varepsilon)$.
    Write $v\vcentcolon =w+\alpha'(0) u = \frac{\dd}{\dd t}\vert_{t=0}\Phi_t$.
    By minimality and since $\Phi_t\vert_{U_1}$ is isometric to $\Phi\vert_{U_1}$, we have
    \begin{align}
          0&=
       \left.\frac{\dd}{\dd t}\right\vert_{t=0} \frac{1}{4}\int_{D\setminus U_1} |\vec H_{\Phi_t}-H_0|^2 \dd\mu_{\Phi_t}  \\
       &=  \int_{D\setminus U_1} \frac{1}{2} (\vec H-H_0\nu)\cdot \Delta v  + \frac{1}{4}\left((3 H - 4H_0)\nabla \nu - \vec H \times \nabla^\perp \nu -(2H_0 H -H_0^2)\nabla \Phi\right)\cdot \nabla v\dd z\,.
    \end{align}
    where we used \Cref{lem:local_first_vari} with $\xi = \chi_{D\setminus {U}_1}$.
    In $U_1$, we have that $v \equiv  a$ is constant, so we may also integrate over all of $D$.
    Hence, by \eqref{eq:def_dE}, we have shown that
    \begin{align}
         0= \delta\calE(\Phi).(w+\alpha'(0) u)\,.
    \end{align}
    Since $\mathcal{A}(\Phi_t)=\mathcal{A}(\Phi)$, we may argue similarly to conclude
    \begin{align}
        0=\left.\frac{\dd}{\dd t}\right\vert_{t=0} \mathcal{A}(\Phi_t) = \delta\mathcal{A}(\Phi).(w+\alpha'(0) u)\,,
    \end{align}
    so that by \Cref{lem:measure-ELE}, we obtain
    \begin{align}
        0= \delta\calE(\Phi).(w+\alpha'(0) u) = -\Lambda \delta\mathcal{A}(\Phi).(w+\alpha'(0) u) +  a\cdot v(b) = | a|^2\,.
    \end{align}
    It follows $ a=0$ and the assertion follows in this case.
    
    If \eqref{eq:constr_nondeg} is satisfied for $i_0\neq i$, we consider $\Psi\vcentcolon = \Phi^{i_0}$ and, with $v$, $\beta$ as in \eqref{eq:def_Psi_t}
    \begin{align}
        \Phi_t \vcentcolon = \Phi+ tw,\quad \Psi_t \vcentcolon = \Psi + \beta(t)v \,.
    \end{align}
    As above, may use minimality, \Cref{lem:local_first_vari}, and $w\equiv a$ in $U_1$ to conclude
    \begin{align}
        0 = \left.\frac{\dd}{\dd t}\right\vert_{t=0} (\calE(\Phi_t) + \calE(\Psi_t)) = \delta\calE(\Phi).w + \beta'(0)\delta\calE(\Psi).v\,.
    \end{align}
    By \Cref{lem:measure-ELE} and since $\frac{\dd}{\dd t}\vert_{t=0}(\mathcal{A}(\Phi_t)+\mathcal{A}(\Psi_t))=0$, we find
    \begin{align}
        0=  \delta\calE(\Phi).w + \beta'(0)\delta\calE(\Psi).v = -\Lambda (\delta \mathcal{A}(\Phi).w + \beta'(0)\delta\mathcal{A}(\Psi). u) +  a\cdot w(0) = | a|^2\,,
    \end{align}
    and again we obtain $a=0$. \qedhere
 \end{proof}

\section{Regularity}

In this section, we examine the local and global regularity of solutions to \eqref{eq:measure-ELE}.

\subsection{The Euler--Lagrange equations as a system of conservation laws}\label{subsec:conservation_laws}

The key observation is that \eqref{eq:measure-ELE} may be transformed into a system of conservation laws with a Jacobian structure, see \Cref{lem:sys_jacobi} below.
For the Willmore functional, this crucial observation has been obtained in \cite{MR2430975}, see also \cite{MR3524220}, whereas for the Helfrich energy we refer to \cites{MR3518329,MR4076069}, and also \cite{scharrer2025energyquantizationconstrainedwillmore}.
Since the impact of a branch point was not taken into account in deriving the corresponding system of conservation laws, we review the calculations for completeness.
This is generically necessary, since $\Phi$ may fail to be immersed at $z=0$.
However, we emphasize that essentially all of the computations in \cite{MR3524220} are still valid in our case, at least in $D\setminus\{0\}$.

As pointed out in \Cref{rem:existence_volume}, the existence theory can easily be extended to also allow for a constraint on the volume or the total mean curvature.
In this case, the nondegeneracy condition \eqref{eq:constr_nondeg} needs to be replaced by an appropriate condition ensuring that the constraints do not degenerate in order to deduce the Euler--Lagrange equations.

For possible future applications, we thus consider \eqref{eq:measure-ELE} with additional Lagrange multipliers corresponding to the volume and the total mean curvature in the sequel. We investigate the local regularity of critical points. Throughout this section, by \Cref{lem:exist_conf_repara}, we  thus assume that we have $\Phi\in \mathcal{F}_D$ conformal, possibly with a single branch point at $z=0$, that satisfies
\begin{align}\label{eq:measure_ELE_1}
    \delta\calE(\Phi) +\Lambda \delta\mathcal{A}(\Phi) + \rho\delta\mathcal{V}(\Phi) +\vartheta \delta\mathcal{T}(\Phi)  = n \circ\Phi \sigma +  a \delta_0 \quad\text{ in }\mathcal{D}'(D;\R^3)\,,
\end{align}
for some $\Lambda,\rho,\vartheta\in\R$. Here $\delta \mathcal{V}(\Phi)$ and $\delta\mathcal{T}(\Phi)$ are the distributions
\begin{align}
    \delta\mathcal{V}(\Phi).w &\vcentcolon= \frac{1}{6} \int_D w\cdot \Div\left[\Phi\times \nabla^\perp\Phi\right]\dd z\, \\
    \delta\mathcal{T}(\Phi).w &\vcentcolon = \int_D (\nabla \nu + H\nabla \Phi)\cdot \nabla w\dd z\qquad \text{ for all }w\in C_c^\infty(D;\R^3)\,.
\end{align}
By \cite{MR4076069}*{(4.6) and (4.9)}, translated to our notation, these correspond to the first variation of total mean curvature and volume, cf.\ \Cref{rem:existence_volume}, for test functions $w$ supported away from the branch point $z=0$.

By \eqref{eq:def_dE} and \eqref{eq:def_dA} and integration by parts in the term with $H_0\nu \cdot\Delta w$, this is equivalent to
\begin{align}
    &\frac{1}{2}\Div\Big[ 2 \nabla \vec H - 3 H \nabla \nu + \vec H \times \nabla^\perp \nu+(2 H_0-4\vartheta)\nabla \nu \\
    &\qquad\quad +\big( (2H_0-4\vartheta) H - H_0^2-4\Lambda\big)\nabla \Phi+ \frac{2\rho}{3} \Phi \times \nabla^\perp\Phi \Big]  = 2n\circ\Phi \dd  \sigma + 2a\delta_0\label{eq:measure-ELE-mal-4}
\end{align}
in $\mathcal{D}'(D;\R^3)$ for some $\Lambda,\rho\in\R$.
To shorten the computations, we write
\begin{align}\label{eq:def_vec_W}
    W \vcentcolon = \frac{1}{2}\left[2\nabla \vec H-3 H\nabla \nu+\vec H\times \nabla^\perp\nu\right]
\end{align}
for the term in \eqref{eq:measure-ELE-mal-4} that arises from the Willmore energy.
First, we exploit the divergence structure of \eqref{eq:measure-ELE-mal-4}. Here and in the following, $L^{2,\infty}(D;\R^3)$ is the weak-$L^2$ Marcinkiewicz space, see \cite{MR2843242}*{Section A.2}. Note that for any $n\in\N$ we have
\begin{align}
    L^2(D;\R^n) \subset L^{2,\infty}(D;\R^n)\subset L^p(D;\R^n)
    \quad\text{ for all }1\leq p<2\,.
\end{align}
\begin{lemma}\label{lem:L_V}
    There exist $ L \in L^{2,\infty}(D;\R^3)$ and $V\in W^{1,p}(D;\R^3)$ for all $p\in [1,2)$ such that in $D$ we have distributionally
    \begin{align}
        \Delta V &= - \frac12 \Div\left[(2H_0-4\vartheta)\nabla \nu + \big( (2H_0-4\vartheta) H-H_0^2-4\Lambda\big)\nabla \Phi +\frac{2\rho}{3}\Phi\times\nabla^\perp\Phi \right]\\
        &\quad + 2n \circ\Phi\dd\sigma + 2 a\delta_0\,, \label{eq:V}\\
        \nabla^\perp  L &=  W- \nabla V\,.
    \end{align}
\end{lemma}
 \begin{proof}
This follows from \eqref{eq:measure-ELE-mal-4}, classical elliptic theory  (taking, for instance, homogeneous Dirichlet boundary conditions), see \cite{MR3156649}*{Theorem 1.2.2}, and the weak Poincar\'e lemma, see \cite{MR2843242}*{Lemma A.1}.
 \end{proof}

\begin{lemma}\label{lem:XY}
    There exist $ X\in W^{2,p}(D;\R^3)$, $Y\in W^{2,p}(D)$ for all $p\in[1,2)$ weakly solving
    \begin{align}
        \Delta X &= \nabla V\times \nabla \Phi\,, \label{eq:X}\\
        \Delta Y &= \nabla V \cdot \nabla \Phi\,, \label{eq:Y}
    \end{align}
    in $D$ and such that we have distributionally
    \begin{align}
        \Div\left[\langle L, \nabla^\perp\Phi\rangle -\nabla Y\right] &=0\,,\label{eq:L_jacobi_system_1}\\
        \Div\left[  L \times \nabla^\perp\Phi + \vec H\times \nabla \Phi - \nabla  X\right]&= 0\,.\label{eq:L_jacobi_system_2}
    \end{align}
\end{lemma}
\begin{proof}
    The existence and regularity of $ X, Y$ is classical elliptic theory.
    For \eqref{eq:L_jacobi_system_1}, first note that we have distributionally
    \begin{align}\label{eq:lem:divL_D_phi}
        \Div\left[ \langle  L,\nabla^\perp\Phi\rangle\right] &=  -  \nabla^\perp  L\cdot \nabla \Phi  = -  W\cdot\nabla \Phi  + \Delta Y\,.
    \end{align}
    Indeed this may be checked by approximating $\Phi$ by smooth functions $\Phi^k\in C_c^\infty(D;\R^3)$ strongly in $W^{2,2}$ and weak-$*$ in $W^{1,\infty}$: For any $ w \in C_c^\infty(D)$, we have
    \begin{align}
       &\Big\langle \langle  L , \nabla^\perp \Phi\rangle, \nabla w\Big\rangle =\int  L \cdot \big(-\partial_2 \Phi \partial_1 w + \partial_1 \Phi \partial_2 w\big)\dd z = \lim_{k\to\infty}\int  L \cdot\big( \partial_2(\partial_1 \Phi^k w) - \partial_1(\partial_2 \Phi^kw)\big)\dd z \\
        &= \lim_{k\to\infty} \Big(\langle \partial_1  L, \partial_2\Phi^k w\rangle - \langle \partial_2  L,\partial_1\Phi^k w\rangle\Big) = \lim_{k\to\infty} \langle \nabla^\perp  L, \nabla \Phi^k w\rangle = \langle  \nabla^\perp  L\cdot\nabla \Phi, w\Big\rangle\,,\label{eq:divLnablaphi_approx}
    \end{align}
    using that $\nabla^\perp L\in W^{-1,2}+L^1$ in the last step.
    The computations in \cite{MR3524220}*{(5.217)--(5.218)} are valid almost everywhere in $D\setminus\{0\}$. Together with the observation that $\langle \nabla H,\nu\rangle\cdot \nabla \Phi=0$ in $\mathcal{D}'(D)$ we have
    \begin{align}
        W\cdot\nabla \Phi=0\quad \text{ in }\mathcal{D}'(D)\,,
    \end{align}
    so \eqref{eq:L_jacobi_system_1} follows from \eqref{eq:lem:divL_D_phi}.
    For \eqref{eq:L_jacobi_system_2}, we note that as in \eqref{eq:divLnablaphi_approx} we find
    \begin{align}
        \Div\left[ L \times \nabla^\perp \Phi \right] =\nabla\Phi\times \nabla^\perp L = \nabla \Phi\times W + \Delta  X\,.
    \end{align}
    As in \cite{MR3524220}*{(5.220)--(5.221)}, we have
    \begin{align}
        \nabla \Phi\times W = \langle \nabla H,\nabla^\perp\Phi\rangle = -  \Div\left[\vec H\times \nabla \Phi\right]\,,
    \end{align}
    where the last equality is a direct computation, using the symmetry of the second fundamental form and the fact that $(\partial_1\Phi,\partial_2 \Phi, \nu)$ is a positively oriented orthonormal basis.
\end{proof}

In the following lemma, $W^{1,(2,\infty)}$ denotes the space of functions in the weak-$L^2$ Marcinkiewicz space $L^{2,\infty}$ with distributional gradient also in $L^{2,\infty}$, see \cite{MR2843242}*{Appendix A.2}.
\begin{lemma}
     There exist $ R \in W^{1,(2,\infty)}(D;\R^3), S\in W^{1,(2,\infty)}(D)$ such that in $D$ we have distributionally
\begin{align}
        \nabla^\perp  R & =  L \times \nabla^\perp\Phi + \vec H\times \nabla \Phi- \nabla  X\,,\label{eq:nabla_perp_R}\\
        \nabla^\perp S & =  \langle  L, \nabla^\perp\Phi\rangle -\nabla Y\,.\label{eq:nabla_perp_S}
        \end{align}
\end{lemma}
\begin{proof}
    This follows from \Cref{lem:XY} and the weak Poincar\'e lemma, see, for instance \cite{MR2843242}*{Lemma A.1}.
\end{proof}
The following structural relation between the potentials is crucial for the regularity analysis.
\begin{lemma}\label{lem:sys_jacobi}
    The triple $( R, S, \Phi)$ is a distributional solution to  the elliptic system
    \begin{align}
        \label{eq:cons_law_R}\Delta R &= \nabla R\times \nabla^\perp \nu - \langle\nabla S,\nabla^\perp \nu\rangle + \Div\left[\langle \nabla Y,\nu\rangle + \nu \times \nabla  X\right]\,,\\
        \label{eq:cons_law_S}\Delta S &= \nabla  R\cdot\nabla^\perp \nu-\Div\left[\langle \nu,\nabla X\rangle\right]\,,\\
        \label{eq:cons_law_Phi}\Delta \Phi &=
        \frac{1}{2}\left(\nabla\Phi \times\nabla^\perp R + \langle \nabla S,\nabla^\perp\Phi\rangle - \langle \nabla^\perp Y,\nabla^\perp\Phi\rangle + \nabla \Phi\times \nabla  X\right)\,.
    \end{align}
    Moreover, we have
    \begin{align}\label{eq:laplace_normal}
    	\Delta \nu &= \nabla^{\perp}\nu \times \nabla \nu - \Div\left[ \vec H\times \nabla^\perp\Phi\right]\,.
    \end{align}
\end{lemma}
\begin{proof}
    From \eqref{eq:nabla_perp_R}, we first obtain
    \begin{align}
        \nabla R =  L \times \nabla\Phi-\vec H\times\nabla^\perp\Phi + \nabla^\perp  X\,.\label{eq:nabla_R_1}
    \end{align}
    On the other hand, as in \cite{MR3524220}*{(5.224)-(5.227)}, we compute
    \begin{align}
        \nu \times \nabla^\perp  R&= \nu \times ( L \times \nabla^\perp\Phi)+ \nu \times (\vec H\times \nabla \Phi) - \nu \times \nabla  X \\
        &=-\langle \nu,\langle  L,\nabla^\perp\Phi\rangle\rangle +  L\times \nabla \Phi - \vec H\times \nabla^\perp\Phi - \nu \times \nabla  X \\
        &= - \langle \nabla^\perp S +\nabla Y ,\nu\rangle +  L\times \nabla \Phi - \vec H\times \nabla^\perp\Phi - \nu \times \nabla  X \\
        &= \nabla  R -\nabla^\perp  X - \langle \nabla^\perp S +\nabla Y ,\nu\rangle  -\nu \times \nabla  X\,.\label{eq:nabla_R_3}
    \end{align}
   Equation \eqref{eq:cons_law_R} follows upon taking the divergence.
   Now, from \eqref{eq:nabla_R_3}, we conclude
    \begin{align}
        \langle \nu,\nabla^\perp  R\rangle &= - \langle\nu,\nu \times \nabla R\rangle + \langle \nu, \langle -\nabla S+\nabla^\perp Y, \nu\rangle\rangle + \langle \nu, \nu \times\nabla^\perp X\rangle - \langle \nu, \nabla  X\rangle \\
        &= -\nabla S + \nabla^\perp Y - \langle \nu,\nabla X\rangle\,.
    \end{align}
    Again, \eqref{eq:cons_law_S} follows from taking the divergence.
    For \eqref{eq:cons_law_Phi}, we note that as in \cite{MR3524220}*{(5.229)}, using \eqref{eq:nabla_perp_S}, we have
    \begin{align}
        \nabla \Phi\times\nabla^\perp R&= \nabla\Phi\times( L \times\nabla^\perp\Phi) - H\nabla\Phi \times\nabla^\perp\Phi-\nabla\Phi\times \nabla  X\\
        &= - \big\langle \langle L, \nabla \Phi\rangle,\nabla^\perp\Phi\big\rangle
        + 2\Delta \Phi  - \nabla\Phi\times\nabla X \\
        &= - \langle \nabla S - \nabla^\perp Y,\nabla^\perp\Phi\rangle
        +2\Delta\Phi - \nabla\Phi\times\nabla X\,.
    \end{align}
    Equation \eqref{eq:laplace_normal} is exactly \cite{MR4706029}*{(A.3)}, using that in $\R^3$ the interior multiplication of $2$-vectors used in \cite{MR4706029} simplifies to the cross product, see also \cite{MR3518329}*{Remark 2.1}.
    The factor 1 instead of 2 is due to our different convention for the mean curvature.
\end{proof}

\subsection{Local regularity of critical points}

We can now discuss the regularity of branched solutions to \eqref{eq:measure-ELE-mal-4}.

\begin{theorem}\label{thm:reg_meas_ELE-2}
Suppose that $\Phi\in \mathcal{F}_D$ is conformal, possibly with a single branch point at $z=0$, and solves \eqref{eq:measure-ELE-mal-4}
    in $D$.
    Then we have $\Phi \in W^{3,p}_{\mathrm{loc}}(B_{1/4}(0);\R^3),$ $\nu \in W^{2,p}_{\mathrm{loc}}(B_{1/4}(0);\R^3)$ for all $p\in[1,2)$.
    In particular, $\Phi$ is $C^{1,\alpha}$ and $\nu$ is $C^\alpha$ throughout the possible branch point $z=0$ for all $\alpha<1$. 
    Moreover, if $\sigma=0$, $ a=0$ in \eqref{eq:measure-ELE-mal-4}, and $\Phi$ is not branched, then $\Phi$ is smooth in $B_{1/4}(0)$.
\end{theorem}
Passing to a smaller disk $B_{1/4}(0)$ in \Cref{thm:reg_meas_ELE-2} is done for technical convenience and the statement remains valid with $B_{1/4}(0)$ replaced by the full disk $D$. Since our interest lies solely in local regularity, we do not pursue this generalization.

Away from the branch points, the conformal factor is bounded from below and the improved regularity yields that we may write the Euler--Lagrange equations \eqref{eq:measure-ELE-mal-4} in its classical strong  scalar form.

\begin{cor}
    In the setting of \Cref{thm:reg_meas_ELE-2}, if there are no branch points in $D$ and after potentially reversing the orientation, we have
    \begin{align}
        \Delta_g H + \big(\frac12 H^2- 2K\big) H + \big(2 H_0-4\vartheta\big) K -\big(2\Lambda+\frac12 H_0^2\big) H - 2\rho =  2e^{-2\lambda}\sigma \qquad \text{ in }\mathcal{D}'(B_{1/4}(0))\,,
    \end{align}
    where $\Delta_g$ is the Laplace--Beltrami operator.
\end{cor}

\begin{proof}[Proof of \Cref{thm:reg_meas_ELE-2}]
    We first discuss some initial regularity. By \Cref{lem:XY}, elliptic theory, and interior regularity,
    there exists a unique solution $ R_0\in W^{2,p}_{\mathrm{loc}}\cap W^{1,p}(D;\R^3)$ and $S_0\in W^{2,p}_{\mathrm{loc}}\cap W^{1,p}(D)$ for all $p\in[1,2)$ to
    \begin{align}\label{eq:system_R0}
     &\begin{cases}
         \Delta  R_0  = \Div \Bigl[\langle \nabla Y,\nu\rangle + \nu\times \nabla X\Bigr] & \text{in } D \\
    		\phantom{\Delta} R_0  = R & \text{on } \partial D\,,
     \end{cases}\\
     \end{align}
     and
     \begin{align}
     \begin{cases}
            \Delta  S_0  = -\Div\Bigl[\langle\nu, \nabla X\rangle\Bigr] & \text{in } D \\
    		\phantom{\Delta} S_0  = S & \text{on } \partial D\,.
     \end{cases}\label{eq:system_S0}
	\end{align}
     Since $\nabla  R\in L^{2,\infty}(D;(\R^3)^2)$, we conclude $S-S_0\in W^{1,2}(D)$ by \eqref{eq:cons_law_S} and integrability by compensation, see \cite{MR3524220}*{Theorem 2.28}.
     Similarly, \eqref{eq:cons_law_R} implies $ R- R_0\in W^{1,2}(D;\R^3)$, so that since
    \begin{align}
    \nabla  R_0\in L^q_{\mathrm{loc}}(D;(\R^3)^2), \nabla S_0\in L^q_{\mathrm{loc}}(D;\R^2) \qquad \text{for all } q\in[1,\infty)\label{eq:nabla_R0}
    \end{align}
    by Sobolev embedding, we have $ R \in W^{1,2}_{\mathrm{loc}}(D;\R^3)$ and $S\in W^{1,2}_{\mathrm{loc}}(D)$.

    \begin{step}[Higher integrability for $\nabla  R, \nabla S$]\label{step:reg-1}
        We claim that there exists $p_0>2$ such that
    \begin{align}\label{eq:RS_higher_reg}
		\nabla  R \in L_{\mathrm{loc}}^{p_0}(B_{1/2}(0); (\mathbb R^3)^2)\,, \qquad \nabla S \in L_{\mathrm{loc}}^{p_0}(B_{1/2}(0); \mathbb R^2)\,.
	\end{align}
     To prove \eqref{eq:RS_higher_reg}, we first show that there exist $r_0\in (0,1/4)$ such that
	\begin{align} \label{eq:morrey}
	\sup_{r<r_0, \, z_0 \in B_{1/2}(0)} r^{-1} \int_{B_r(z_0)} \Big(|\nabla  R|^2 + |\nabla S|^2\Big) \dd z < \infty\,.
	\end{align}
    Let $\varepsilon_0 > 0$ to be chosen in the sequel.
    There exists $0 < r_0 < 1/4$ such that
	\begin{align} \label{eq:small_nabla_n}
		\sup_{z_0 \in B_{1/2}(0)} \int_{B_{r_0}(z_0)} |\nabla \nu|^2 \, \dd z < \varepsilon_0\,.
	\end{align}
    Fix $z_0\in B_{1/2}(0)$ and $0<r\leq r_0$.
    Using \eqref{eq:nabla_R0} with $q=4$, the Cauchy--Schwarz inequality implies that
	\begin{align} \label{eq:nablaR0_estimate}
		&\int_{B_r(z_0)}\Big(|\nabla  R_0|^2 + |\nabla S_0|^2\Big) \,\dd z \\
        &\leq \sqrt{4\pi r^2} \Biggl(\Bigl(\int_{B_1(0)} |\nabla  R_0|^4 \, \dd z\Bigl)^{1/2} +\Bigl(\int_{B_1(0)} |\nabla  S_0|^4 \, \dd z\Bigl)^{1/2} \Biggr) =\vcentcolon r C_1\,.
	\end{align}
    Note that $C_1$ does not depend on $z_0$ and $r$.
	We now let $\Psi_{ R}\colon B_r(z_0)\to\R^3$ and $\Psi_S\colon B_r(z_0)\to\R$ be the solutions of
 	\begin{align}
    \label{eq:system_PsiR}
    \begin{cases}
        \Delta \Psi_{ R} = \langle \nabla^\bot \nu, \nabla S \rangle + \nabla^\bot \nu \times \nabla  R & \text{in } B_{r}(z_0) \\
		\phantom{\Delta}\Psi_{ R}  = 0 & \text{on } \partial B_{r}(z_0)
    \end{cases}
	\end{align}
	and
	\begin{align}
    \label{eq:system_PsiS}
    \begin{cases}
        \Delta \Psi_{S}  = \langle\nabla^\bot \nu,\nabla  R\rangle & \text{in } B_{r}(z_0) \\
		\phantom{\Delta}\Psi_{S} = 0 & \text{on } \partial B_{r}(z_0)\,.
    \end{cases}
    \end{align}
    Again, integrability by compensation implies $\Psi_{ R} \in W^{1,2}(B_r(z_0);\R^3)$ and $\Psi_S\in W^{1,2}(B_r(z_0))$.
	Then, by \eqref{eq:cons_law_R} and \eqref{eq:cons_law_S}, the harmonic functions
    \begin{align}\label{eq:def_harmonic}
    v_{ R} \vcentcolon =  R -  R_0 - \Psi_{ R}\,,\qquad v_{S} \vcentcolon = S - S_0- \Psi_S
    \end{align}
    minimize the Dirichlet energy in the class of all $W^{1,2}$-functions with boundary values
	\begin{align}
	v_{ R} =  R -  R_0\,, \quad v_{S} = S-S_0 \qquad \text{on } \partial B_r(z_0)\,.
	\end{align}
    Take $0<\rho\leq r<r_0$.
    By monotonicity (see \cite{MR3524220}*{Lemma 5.65}), and minimality
	\begin{align}
		\rho^{-2}\int_{B_{\rho}(z_0)} \Big(|\nabla v_{ R}|^2 + |\nabla v_{S}|^2\Big) \dd z & \leq r^{-2}\int_{B_{r}(z_0)} \Big(|\nabla ( R -  R_0)|^2 + |\nabla (S-S_0)|^2\Big) \dd z \\
		& \leq 2r^{-2} \int_{B_{r}(z_0)} \Big(|\nabla  R |^2 + |\nabla S|^2 \Big) \dd z + 2r^{-1} C_1\,,
    \label{eq:harmonic_rest_estimate}
	\end{align}
    using \eqref{eq:nablaR0_estimate} in the last step.
    Using first
	Wente's theorem (see for instance \cite{MR3524220}*{Theorem 2.25}) for \eqref{eq:system_PsiR} and \eqref{eq:system_PsiS}, and then \eqref{eq:small_nabla_n}, we find
	\begin{align}
		\int_{B_{r}(z_0)} \Big(|\nabla \Psi_{ R}|^2 + |\nabla \Psi_S|^2 \Big)\dd z  &\leq C_2\int_{B_{r_0}(z_0)} |\nabla \nu|^2 \dd z \int_{B_{r}(z_0)} \Big(|\nabla  R|^2 + |\nabla S|^2 \Big) \dd z \\
		&  \leq C_2\varepsilon_0\int_{B_{r}(z_0)} \Big(|\nabla  R|^2 + |\nabla S|^2 \Big)\dd z \label{eq:Wente_estimate}
	\end{align}
	for some universal constant $C_2 \in (0,\infty)$.
     Recalling \eqref{eq:def_harmonic} and combining \eqref{eq:nablaR0_estimate}, \eqref{eq:harmonic_rest_estimate}, and \eqref{eq:Wente_estimate}, we find
	\begin{align}
		& \int_{B_{\rho}(z_0)} \Big(|\nabla  R|^2 + |\nabla S|^2 \Big) \dd z \\
		& \quad \leq 3\int_{B_{\rho}(z_0)}\Big(|\nabla \Psi_{ R}|^2 + |\nabla \Psi_S|^2 \Big)\dd z + 3\int_{B_{\rho}(z_0)} \Big(|\nabla v_{ R}|^2 + |\nabla v_{S}|^2 \Big) \dd z \\
        &\quad \qquad+ 3\int_{B_{\rho}(z_0)}\Big(|\nabla  R_0|^2 +|\nabla S_0|^2\Big)\dd z \\
		& \quad \leq \Big(3C_2\varepsilon_0 + \frac{6\rho^2}{r^2}\Big) \int_{B_{r}(z_0)} \Big( |\nabla  R|^2 + |\nabla S|^2 \Big) \dd z + 9r C_1\,.
	\end{align}
    Now \eqref{eq:morrey} follows from a classical iteration lemma, see for instance \cite{MR3099262}*{Lemma 5.13}, if we take $\varepsilon_0 = \varepsilon_0(C_2)>0$ small enough.
    
	  Using \eqref{eq:cons_law_R}, \eqref{eq:cons_law_S}, \eqref{eq:system_R0}, \eqref{eq:system_S0}, and Cauchy--Schwarz\,, \eqref{eq:morrey} implies
	\begin{align}
		\sup_{r<r_0, z_0 \in B_{1/2}(0)} r^{-1/2} \int_{B_r(z_0)} \Big(|\Delta( R -  R_0)| + |\Delta (S-S_0)|\Big) \dd z < \infty\,.
	\end{align}
	Then a classical estimate on Riesz potentials \cite{MR458158} gives
	\begin{align}
		\nabla ( R -  R_0) \in L_{\mathrm{loc}}^{p_0}(B_{1/2}(0); (\mathbb R^3)^2)\,, \qquad \nabla (S-S_0) \in L_{\mathrm{loc}}^{p_0}(B_{1/2}(0); \mathbb R^2)\,,
	\end{align}
	for all $p_0\in (2,3)$,  see also \cite{MR4706029}*{Proposition A.1} for the detailed argument.
	From \eqref{eq:nabla_R0}, we conclude \eqref{eq:RS_higher_reg}. $\hfill\diamond$
    \end{step}
    
    \begin{step}[Higher integrability of $\nabla \nu$]\label{step:reg-2}
        We now show that there exist $p_1 \in (2,p_0]$ with
    \begin{align}\label{eq:higher_int_n}
        \nabla \nu \in L^{p_1}_{\mathrm{loc}}(B_{1/4}(0);(\R^3)^2)\,.
    \end{align}
    To that end, we first note that using $e^{2\lambda}\vec H = \Delta\Phi$ in $D\setminus \{0\}$ and \eqref{eq:cons_law_Phi} yields
    \begin{align}
        2 e^{2\lambda}\vec H \times \nabla^\perp\Phi = \Big(\nabla \Phi\times \nabla^\perp R + \langle \nabla S, \nabla^\perp\Phi\rangle - \langle \nabla^\perp Y,\nabla^\perp\Phi\rangle + \nabla \Phi\times \nabla X\Big)\times \nabla^\perp\Phi\,,\label{eq:H_x_Dphi}
    \end{align}
    thus, using $|\nabla \Phi|^2=e^{2\lambda}$, for a.e.\ $z\in D$ we may estimate
    \begin{align}\label{eq:3.42a}
        2|\vec H \times \nabla^\perp\Phi|\leq |\nabla^\perp R| + |\nabla S| + |\nabla^\perp Y|+|\nabla X|\,.
    \end{align}
    By \Cref{step:reg-1} and \Cref{lem:XY}, this implies that $\vec H\times \nabla^\perp\Phi\in L^{p_0}_{\mathrm{loc}}(B_{1/2}(0);(\R^3)^2)$.
    Hence, by elliptic regularity, the solution $\nu_0$ of
    \begin{align}\label{eq:system_n0}
		\begin{cases}
			\Delta \nu_0 =-\Div\left[ \vec H\times \nabla^\perp\Phi\right] & \text{in }B_{1/2}(0)\\
			\phantom{\Delta}\nu_0 = 0 & \text{on }\partial B_{1/2}(0)
		\end{cases}
	\end{align}
    satisfies $\nu_0\in W^{1,p_0}_{\mathrm{loc}}(B_{1/2}(0);\R^3)$ with $p_0>2$ as in \Cref{step:reg-1}.
    Hölder's inequality implies that, after possibly reducing $r_0\in (0,1/4)$, there exist $\alpha_0 >0$ and $C_3\in(0,\infty)$ such that we have
    \begin{align}\label{eq:nabla_n0_estimate}
        \int_{B_r(z_0)} |\nabla \nu_0|^2\dd z\leq C_3 r^{\alpha_0}\,,
    \end{align}
    for all $z_0\in B_{1/4}(0)$ and $r\in (0,r_0]$.
    Fix some $z_0\in B_{1/4}(0)$, $r\in (0,r_0]$, and let $\Psi_{\nu}$ be the solution of
    \begin{align}
    \label{eq:system_Psi_n}
    \begin{cases}
        \Delta \Psi_{\nu}  = \nabla^\perp \nu \times \nabla \nu & \text{in } B_{r}(z_0) \\
		\phantom{\Delta}\Psi_{\nu} = 0 & \text{on } \partial B_{r}(z_0)\,.
    \end{cases}
    \end{align}
    Integrability by compensation yields $\Psi_{\nu}\in W^{1,2}(B_r(z_0);\R^3)$.
    By \eqref{eq:laplace_normal}, the function $v_{\nu} \vcentcolon = \nu - \nu_0 - \Psi_{\nu}$ is harmonic with boundary value $v_{\nu} = \nu - \nu_0$ on $\partial B_r(z_0)$.
    As in \eqref{eq:harmonic_rest_estimate}, for $0<\rho\leq r<r_0$, using \eqref{eq:nabla_n0_estimate} we find
    \begin{align}\label{eq:harmonic_rest_estimate_2}
         \rho^{-2} \int_{B_{\rho}(z_0)} |\nabla v_{\nu}|^2\dd z \leq r^{-2} \int_{B_r(z_0)} |\nabla(\nu - \nu_0)|^2\dd z \leq 2 r^{-2} \int_{B_r(z_0)} |\nabla \nu|^2\dd z + 2 r^{\alpha_0-2}C_3\,.
    \end{align}
    By Wente's theorem (see \cite{MR3524220}*{Theorem 3.7}) with $C_2$ as in \eqref{eq:Wente_estimate}, we have
	\begin{align}
		\int_{B_{r}(z_0)} |\nabla \Psi_{\nu}|^2 \dd z  &\leq C_2\int_{B_{r_0}(z_0)} |\nabla \nu|^2 \dd z \int_{B_{r}(z_0)} |\nabla \nu|^2 \dd z \leq C_2\varepsilon_0\int_{B_{r}(z_0)}|\nabla \nu|^2 \dd z\,,
	\end{align}
    using \eqref{eq:small_nabla_n} in the last step.
    Combining this with \eqref{eq:nabla_n0_estimate} and \eqref{eq:harmonic_rest_estimate_2}, we have
	\begin{align}
		\int_{B_{\rho}(z_0)} |\nabla \nu|^2 \dd z &\leq 3\int_{B_{\rho}(z_0)} |\nabla \Psi_{\nu}|^2\dd z + 3\int_{B_{\rho}(z_0)}|\nabla v_{\nu}|^2\dd z + 3\int_{B_{\rho}(z_0)} |\nabla \nu_0|^2\dd z \\
		& \leq  \Big(3 C_2 \varepsilon_0 + \frac{6\rho^2}{r^2}\Big) \int_{B_r(z_0)}|\nabla \nu|^2\dd z + 9C_3 r^{\alpha_0}\,.
	\end{align}
    As before, \cite{MR3099262}*{Lemma 5.13} implies
    that for some $\alpha_1>0$ we have
    \begin{align}
        \sup_{r<r_0,\, z_0\in B_{1/4}(0)} r^{-\alpha_1} \int_{B_r(z_0)} |\nabla\nu|^2\dd z <\infty\,,
    \end{align}
    if we take $\varepsilon_0(C_2,\alpha_0)>0$ in \eqref{eq:small_nabla_n} sufficiently small.
    Now, \eqref{eq:higher_int_n} follows as in \Cref{step:reg-2}, using \eqref{eq:laplace_normal} and that $\nu_0\in W^{1,p_0}_{\mathrm{loc}}(B_{1/4}(0);\R^3)$. $\hfill\diamond$
    \end{step}
    
    \begin{step}[Bootstrap]
        We may now complete the proof.
    By Steps \ref{step:reg-1} and \ref{step:reg-2}, there exist $q_0\in (2,\infty)$ (in fact $q_0\vcentcolon = p_1$ works) with
    \begin{align}\label{eq:boostrap_0}
        \nabla  R, \nabla \nu\in L^{q_0}_{\mathrm{loc}}(B_{1/4}(0);(\R^3)^2)\,, \quad \nabla S\in L^{q_0}_{\mathrm{loc}}(B_{1/4}(0);\R^2)\,.
    \end{align}
    Using Hölder's inequality, \eqref{eq:cons_law_R}, \eqref{eq:cons_law_S}, \eqref{eq:laplace_normal}, \eqref{eq:system_R0}, \eqref{eq:system_S0}, and \eqref{eq:system_n0}, we have
    \begin{align}
        \Delta( R- R_0), \Delta(\nu - \nu_0)\in L^{q_0/2}_{\mathrm{loc}}(B_{1/4}(0);\R^3)\,, \quad \Delta (S-S_0)\in L^{q_0/2}_{\mathrm{loc}}(B_{1/4}(0))\,.
    \end{align}
    Elliptic regularity yields
    \begin{align}
         R- R_0, \nu - \nu_0 \in W^{2,q_0/2}_{\mathrm{loc}}(B_{1/4}(0);\R^3)\,, \quad S-S_0\in W^{2,q_0/2}_{\mathrm{loc}}(B_{1/4}(0))\,.
    \end{align}
    If $q_0/2<2$, by Sobolev embedding
    \begin{align}\label{eq:bootstrap_1}
        \nabla( R- R_0), \nabla(\nu - \nu_0)\in L^{q_1}_{\mathrm{loc}}(B_{1/4}(0);(\R^3)^2)\,,\quad \nabla (S-S_0)\in L^{q_1}_{\mathrm{loc}}(B_{1/4}(0);\R^2)\,,
    \end{align}
    with $q_1\vcentcolon = (q_0/2)^* = 2q_0/(4-q_0)>q_0$.
    If $q_0/2\geq 2$, then  \eqref{eq:bootstrap_1} is true for $q_1\vcentcolon =q_0+1$.
    We obtain step by step
    \begin{align}
        &\nabla  R \in L^{q_1}_{\mathrm{loc}}(B_{1/4}(0);(\R^3)^2), \nabla S\in L^{q_1}_{\mathrm{loc}}(B_{1/4}(0);\R^2)\,, && \text{by  \eqref{eq:nabla_R0} and \eqref{eq:bootstrap_1},}\\
        &\nabla \nu_0 \in L^{q_1}_{\mathrm{loc}}(B_{1/4}(0);(\R^3)^2)\,, && \text{by \eqref{eq:3.42a} and \eqref{eq:system_n0},}\\
        &\nabla \nu \in L^{q_1}_{\mathrm{loc}}(B_{1/4}(0);(\R^3)^2)\,, && \text{by \eqref{eq:bootstrap_1}.}
    \end{align}
    Constructing $q_2$ from $q_1$ as above, we obtain a recursive sequence $(q_n)_{n\in\N}$, which is easily shown to be strictly increasing and unbounded, so \eqref{eq:boostrap_0} is true for any $q_0\in (2,\infty)$.
    Recalling $Y\in W^{2,p}(D)$ for all $p\in [1,2)$, from \eqref{eq:cons_law_R} and \eqref{eq:cons_law_S}, we conclude
    \begin{align}
         R \in W^{2,p}_{\mathrm{loc}}(B_{1/4}(0);\R^3)\,, \quad S\in W^{2,p}_{\mathrm{loc}}(B_{1/4}(0)) \text{ for all } p\in [1,2)\,.
    \end{align}
    Now, from \eqref{eq:cons_law_Phi}, \eqref{eq:3.42a}, and \eqref{eq:laplace_normal}, it follows that for all $p\in[1,2)$ we have
    \begin{align}
        \Phi \in W^{3,p}_{\mathrm{loc}}(B_{1/4}(0);\R^3)\,, \quad \nu \in W^{2,p}_{\mathrm{loc}}(B_{1/4}(0);\R^3)\,.
    \end{align}

    For the last part of the statement, note that if $\Phi$ is unbranched, $|\nabla\Phi|^2=e^{2\lambda}$ is bounded from above and below. It thus follows that $H\in L^p(B_{1/4}(0))$ for all $p<\infty$. Consequently, by \eqref{eq:V}, we have $V\in W^{1,p}(B_{1/4}(0);\R^3)$ for all $p<\infty$, so \eqref{eq:X}, \eqref{eq:Y} imply $ X \in W^{2,p}(B_{1/4}(0);\R^3)$ and $Y\in W^{2,p}(B_{1/4}(0))$ for all $p<\infty$.
    Bootstrapping \eqref{eq:cons_law_R}--\eqref{eq:cons_law_Phi} together with \eqref{eq:V}, \eqref{eq:X}, and \eqref{eq:Y} yields the claim.\qedhere
    \end{step}
\end{proof}

\begin{remark}
    In the setting of \Cref{thm:reg_meas_ELE-2}, from \eqref{eq:cons_law_Phi} and since $|\nabla\Phi|^2=e^{2\lambda}$, we conclude
    \begin{align}
        e^{2\lambda}|\vec H| = |\Delta\Phi| \leq \frac{e^\lambda}{2} \big(|\nabla R|+|\nabla S|+|\nabla  X|+|\nabla Y|\big) \quad \text{a.e.\ in }D\,,
    \end{align}
    so that even in the branched case $e^{2\lambda} \vec H\in L^p_{\mathrm{loc}}(B_{1/4}(0);\R^3)$ for all $p<\infty$ follows.
\end{remark}

\subsection{Global regularity of minimizers}
The global regularity of minimizers can be summarized as follows.
\begin{theorem}\label{thm:reg_global}
    Let $ T_* = (f,\Phi^1,\dots,\Phi^N)$ be a minimizer of \Cref{thm:existence_bubble}.
    Then at least  one of the following two alternatives is true.
    \begin{enumerate}
        \item\label{item:thm:reg:1} For all $i=1,\dots,N$, we have $\vec H_{\Phi^i}\equiv 0$ a.e.\ in $\S^2\setminus (\Phi^i)^{-1}(\partial\Omega)$.
        In this case, $\vec H_{\Phi^i}\in L^\infty(\S^2;\R^3)$ and $\Phi^i\in W^{2,p}(\S^2;\R^3)$ for all $p\in[1,\infty)$, $i=1,\dots,N$.
        \item\label{item:thm:reg:2} For all $i=1,\dots,N$, we have the following.
        \begin{enumerate}
            \item\label{item:thm:reg:2a}
            There exists $\Lambda^i\in\R$, a positive Radon measure $\sigma^i$ supported in $(\Phi^i)^{-1}(\partial\Omega)$, and $ a^{i,1},\dots, a^{i,N^i}\in\R^3$ with $ a^{i,k} = 0$ if $\Phi^i(b^{i,k})\in \Omega$ such that
        \begin{align}\label{eq:thm:reg:ELE}
             \delta \calE(\Phi^i) + \Lambda^i \delta\mathcal{A}(\Phi^i) = n\circ\Phi^i \sigma^i + \sum_{k=1}^{N^i}  a^{i,k} \delta_{b^{i,k}}\quad\text{ in $\mathcal{D}'(\S^2;\R^3)$}\,.
        \end{align}
        \item\label{item:thm:reg:2b} We have $\Phi^i\in W^{3,p}(\S^2;\R^3)$, $\nu_{\Phi^i}\in W^{2,p}(\S^2;\R^3)$ for each $p\in [1,2)$.
        \item\label{item:thm:reg:2c} $\Phi^i$ is smooth in $(\Phi^i)^{-1}(\Omega) \setminus \{b^{1,1}, \dots, b^{i,N^i}\}$.
        \end{enumerate}
    \end{enumerate}
\end{theorem}

\begin{proof}
    If \eqref{eq:constr_nondeg} is not satisfied, then $\vec H_{\Phi}\equiv 0$ a.e.\ in $\Phi^{-1}(\partial \Omega)$.
    By \cite{MR2472179}*{Theorem 4.1}, for almost every $x\in \S^2$, we have
    \begin{align}
        \vec H_{\Phi}(x) = \begin{cases}
            \vec H_{\partial\Omega}(\Phi(x)) & \text{if } \Phi(x)\in \partial\Omega\,,\\
            0 &\text{otherwise.}
        \end{cases}
    \end{align}
    In particular, it follows that $\vec H_{\Phi}\in L^\infty(D;\R^3)$, so that $\Delta\Phi = e^{2\lambda}\vec H\in L^p(D;\R^3)$ for all $p<\infty$.
    Elliptic regularity yields alternative \eqref{item:thm:reg:1}.

    On the other hand, suppose that \eqref{eq:constr_nondeg} is satisfied and fix $i\in \{1,\dots,N\}$.
    By \Cref{lem:measure-ELE}, we have the weak Euler--Lagrange equations \eqref{eq:measure-ELE} for some $\Lambda^i\in\R$, a Radon measure $\sigma^i$, and $a^{i,1},\dots,a^{i, N^i}\in \R^3$ as statement \eqref{item:thm:reg:2a}, where the fact that $ a^{i,k}=0$ if $\Phi^i (b^{i,k})\in\Omega$ follows from \Cref{lem:residue}.
    The regularity statement \eqref{item:thm:reg:2b} follows from \Cref{thm:reg_meas_ELE-2}, since in local conformal coordinates for $\Phi^i$, \eqref{eq:thm:reg:ELE} implies that \eqref{eq:measure_ELE_1} is satisfied.
    Lastly, in $(\Phi^i)^{-1}(\Omega)\setminus \{b^{i,1},\dots,b^{i,N^i}\}$, $\Phi^i$ is a Lipschitz immersion, the measure terms in \eqref{eq:thm:reg:ELE} vanish, and \eqref{item:thm:reg:2c} follows again from \Cref{thm:reg_meas_ELE-2} after taking local conformal coordinates.
\end{proof}

In terms of equation \eqref{eq:measure-ELE-mal-4}, the regularity in \Cref{thm:reg_meas_ELE-2} is optimal.
\begin{lemma}[Regularity of the inverted catenoid]\label{lem:inv_cat}
    There exists $\Phi\colon D\to\R^3$ a conformal parametrization of (one side of) the inverted catenoid around the origin with the following properties.
    \begin{enumerate}
        \item $\Phi$ solves \eqref{eq:measure-ELE-mal-4} for $H_0=0$ with $\sigma=0$, $ a\neq 0$,
        \item $\Phi \in W^{3,p}(D;\R^3)$ for all $p<2$ and $\Phi\not\in W^{3,2}(D;\R^3)$.
    \end{enumerate}
\end{lemma}
\begin{proof}
    As pointed out in \cite{MR3096502}*{Remark 1.4}, the parametrization $f\colon\R^2\to\R^3$ of the catenoid given in polar coordinates $r>0, \varphi\in [0,2\pi)$ by
    \begin{align}
        f(r,\varphi) = \left( (r+r^{-1}) \cos\varphi, (r+r^{-1}) \sin\varphi, -2 \log r\right)
    \end{align}
    is conformal.
    Inverting in the unit sphere, we obtain a conformal Lipschitz immersion $\Phi (z) \vcentcolon = |f(z)|^{-2} f(z)$, $z\in D$, in particular we have $\Phi\in \mathcal{F}_D$.
    As observed in \cite{MR3096502}*{Remark 1.4}, $\Phi$ solves the Willmore equation in $D$ with some some nonzero Dirac contribution at the origin.
    Hence $\Phi$ solves \eqref{eq:measure-ELE-mal-4} for $H_0=0$, $\sigma=0$, and $ a\neq 0$.
    A short explicit computation yields
    \begin{align}
        \Phi(z) = \left( \frac{z_1(1+r^2)}{(1+r^2)^2 + r^2 \log(r^2)^2}, \frac{z_2(1+r^2)}{(1+r^2)^2 + r^2 \log(r^2)^2}, \frac{-r^2 \log(r^2)}{(1+r^2)^2+r^2 \log(r^2)^2}\right)
    \end{align}
    where $r=|z|$.
    It now suffices to show that $\Phi\not\in W^{3,2}(D;\R^3)$ since $\Phi\in W^{3,p}(D;\R^3)$ for all $p<2$ follows from \Cref{thm:reg_meas_ELE-2}.
    We now write $g(r)$ for the Laplacian of the third component of $\Phi$, i.e.,
    \begin{align}
        g(r) = \Big(\partial_r^2 + r^{-1} \partial_r\Big)~\frac{-r^2 \log(r^2)}{(1+r^2)^2+r^2 \log(r^2)^2}\,.
    \end{align}
    It is a straightforward but tedious computation to check that
    \begin{align}
        \lim_{r\to 0+} g'(r) r = -8\,.
    \end{align}
    It follows that $|g'(r)|\geq c/r$ for $r>0$ small, so that
    \begin{align}\label{eq:nabla_Delta_polar}
        \int_D |\nabla \Delta \Phi_3|^2\dd z =  2\pi \int_0^1 |g'(r)|^2 r\dd r =\infty\,. &\qedhere
    \end{align}
\end{proof}

\subsection{Further applications}\label{sec:futher_applications}
Since our regularity analysis (throughout the branch point) uses only the weak Euler--Lagrange equations \eqref{eq:measure-ELE-mal-4} in local conformal coordinates, \Cref{thm:reg_meas_ELE-2} can be used to prove regularity throughout the branch points in various variational setups without confinement conditions. Recall  that \Cref{thm:reg_meas_ELE-2} yields that any conformal $\Phi\in \mathcal{F}_D$, possibly with a single branch point at $z=0$, that solves \eqref{eq:measure-ELE-mal-4} with $\sigma=0$, i.e.,
\begin{align}
    \delta \E_{H_0}(\Phi) + \Lambda \delta\mathcal{A}(\Phi) + \rho\delta \mathcal{V}(\Phi) + \vartheta \delta \mathcal{T}(\Phi) = a \delta_0 \qquad\text{ in }\mathcal{D}'(D;\R^3)\, \label{eq:branch_constr_helfrich}
\end{align}
satisfies $\Phi\in W^{3,p}(D;\R^3)$ and $\nu\in W^{2,p}(D;\R^3)$ for all $p\in[1,2)$. Applying \Cref{lem:exist_conf_repara}, this arises as the local Euler--Lagrange equations near a possible branch point of the Helfrich energy (or the Willmore energy if $H_0=0$) with possible constraints on the area, volume, or total mean curvature, or any combination of these. In particular, minimizers of $\E_{H_0}$ in $\mathcal{F}_\Sigma$ or in the class of bubble trees with possible constraints on area, volume, or total mean curvature, will satisfy \eqref{eq:branch_constr_helfrich} for some choices of the parameters $H_0,\Lambda, \rho,\vartheta\in\R$, and $a\in \R^3$, provided the constraints do not degenerate. We gather some consequences for the regularity through the branch points of minimizers that were shown to exist in previous works.

Mondino--Scharrer \cite{MR4076069} showed the existence of minimizers for the Helfrich functional in the class of bubble trees with prescribed area and volume, and their smoothness away from the branch points \cite{MR4076069}*{Theorem 1.7}. 
Around a branch point of a bubble $\Phi$ of a minimizer in \cite{MR4076069}*{Theorem 1.7}, by \Cref{lem:exist_conf_repara} we may find local conformal coordinates such that (by \cite{MR4076069}*{Lemma 4.1}), in $\mathcal{D}'(D\setminus\{0\};\R^3)$ Equation \eqref{eq:branch_constr_helfrich} is satisfied (with $\sigma=0$). By the same argument as in \Cref{step:ele-5} in the proof of \Cref{lem:measure-ELE}, we conclude that \eqref{eq:branch_constr_helfrich} is satisfied in $\mathcal{D}'(D;\R^3)$ for $\sigma=0$, $\vartheta=0$, and some $a\in \R^3$. Applying \Cref{thm:reg_meas_ELE-2}, we thus find the following.
\begin{cor}\label{cor:MS_regularity}
	Let $T=(f, \Phi^1, \dots, \Phi^N)$ be a minimizer in \cite{MR4076069}*{Theorem 1.7}. Then for all $i\in \{1,\dots,N\}$, we have $\Phi^i \in W^{3,p}(\S^2;\R^3)$ and $\nu_{\Phi^i} \in W^{2,p}(\S^2;\R^3)$ for all $p\in [1,2)$. In particular, $\Phi^i$ is $C^{1,\alpha}$ and $\nu_{\Phi^i}$ is $C^\alpha$ throughout the branch points.
\end{cor}

In \cite{MR4076072},  Da Lio--Palmurella--Rivi\`ere studied a curvature minimization problem with clamped boundary conditions and proved existence of minimizers (in the class $\mathcal{F}_D$) which are smooth away from the branch points and their $C^{1,\alpha}$-regularity through the branch points. They also showed that, away from the boundary equation \eqref{eq:branch_constr_helfrich} is satisfied with $\sigma=0$, $\rho=\vartheta=0$, $H_0=0$, $a=0$, see \cite{MR4076072}*{Theorem 1.3}. Our analysis provides the optimal interior Sobolev regularity for critical points in their variational setup.
\begin{cor}\label{cor:DLPR_regularity}
	Consider $\Phi\in \mathcal{F}_D$ conformal which is a minimizer in \cite{MR4076072}*{Theorem 1.2}. Then $\Phi\in W^{3,p}_{\mathrm{loc}}(D;\R^3)$ and $\nu \in W^{2,p}_{\mathrm{loc}}(D;\R^3)$ for all $p\in [1,2)$.
\end{cor}

Sequences of solutions to \eqref{eq:branch_constr_helfrich} with $a=0$ have been studied in \cite{scharrer2025energyquantizationconstrainedwillmore}. 
Under suitable assumptions, \cite{scharrer2025energyquantizationconstrainedwillmore}*{Theorem 1.2} yields a limit bubble tree of branched immersions locally solving \eqref{eq:branch_constr_helfrich} in conformal coordinates (possibly with $a\neq 0$). \Cref{thm:reg_meas_ELE-2} gives that any of these immersions is $W^{3,p}$-regular with $W^{2,p}$-regular unit normal for all $p\in [1,2)$.

The highest order terms in \eqref{eq:branch_constr_helfrich} arise from the Willmore energy. For $H_0=\Lambda=\rho=0$, solutions to \eqref{eq:branch_constr_helfrich} are branched Willmore immersions whose regularity has been studied in several works, see \cites{KuSc04,MR2318282} and \cite{MR3096502}. In particular, these works show $C^{1,\alpha}$-regularity for all $\alpha<1$ throughout the branch point which cannot be improved to $C^{1,1}$-regularity due to the inverted catenoid, see \cite{KuSc04}*{p.~337}. \Cref{thm:reg_meas_ELE-2} and \Cref{lem:inv_cat} show the optimal regularity of the inverted catenoid also in the Sobolev scale.
\begin{cor}\label{cor:branch_Willmore_regularity}
	Let $\Phi\in \mathcal{F}_D$ be conformal, possibly with a single branch point at $z=0$, and a Willmore immersion in $D\setminus \{0\}$. Then $\Phi\in W^{3,p}(D;\R^3)$ and $\nu \in W^{2,p}(D;\R^3)$ for all $p\in [1,2)$.
\end{cor}

\section{Rigidity for the confinement problem in the ball}
\subsection{Lower bounds and the confinement problem for arbitrary area}
We consider the minimization of the Helfrich energy under a confinement constraint with the container given by the unit ball:
Let $\calM$ denote the set of smooth, embedded closed surfaces and consider
\begin{align}
    \label{eq:inf_ext} \opt(H_0) := \inf\{ \calE_{H_0}(S) \mid S\in\calM\,, S\subset \overline{B}\}\,.
\end{align}
A first observation is that the energy vanishes if $H_0$ favors small spheres.
More surprisingly, the preference for spheres is rather rigid and extends to the regime $0\leq H_0\leq 2$.

\begin{theorem}\label{thm:main_ext_obstacle}
\begin{enumerate}
    \item For any $H_0\geq 2$ it holds $\opt(H_0)=0$ and the infimum is attained only by spheres $\partial B_r(x)\subset \overline{B}$ of radius $r=\frac{2}{H_0} \leq 1$.
    \item For any $0\leq H_0\leq 2$
    \begin{align}\label{eq:E(B)_minimal}
        \opt(H_0)\geq \calE_{H_0}(\partial B)
        =(2-H_0)^2\pi
    \end{align}
    holds, and the infimum is attained if and only if $S=\partial B$.
    \item For any $H_0<0$ it holds $\opt(H_0)=4\pi$ and no minimizer exists in $\calM$.
\end{enumerate}
\end{theorem}

\begin{proof}
\begin{step}
    The first statement in the theorem follows by a direct computation and the fact that the only closed surfaces with constant mean curvature are spheres.
\end{step}

\begin{step}
    Using the confinement condition $S\subset \overline{B_1(0)}$ and the first variation formula we have
\begin{align}
    2|S| &= \int_S -\vec{H}(x)\cdot x\dd\calH^2(x)\\
    &\leq \int_S |x\cdot\nu(x)|^2 \dd\calH^2(x) + \frac{1}{4}\int_S H^2\dd\calH^2
    \leq |S| + \frac{1}{4}\int_S H^2\dd\calH^2\,.
\end{align}
Since balls minimize the Willmore energy in the class of closed smooth surfaces we deduce
\begin{equation}
    \frac{1}{4}\int_S H^2\dd\calH^2 \geq M:=\max\{4\pi,|S|\}
    \quad\text{ for all }S\in\calM\,,
    \label{eq:PfMainMin-4}
\end{equation}
with equality for $|S|<4\pi$ if and only if $S$ is a ball, and for $|S|\geq 4\pi$ if and only if $S=\partial B$.
\end{step}

\begin{step}
    Even without any assumption on $H_0\in\R$ we obtain
\begin{align}
    &\int_S (H-H_0)^2\dd\calH^2 \nonumber\\
    &\qquad= \big\||H|-|H_0|\big\|_{L^2(S)}^2 +4\int_S (HH_0)_-\dd\calH^2\nonumber\\
    &\qquad\geq \big(\|H\|_{L^2(S)}-\|H_0\|_{L^2(S)}\big)^2+4\int_S (HH_0)_-\dd\calH^2\nonumber\\
    &\qquad= \big(2\sqrt{M}-\|H_0\|_{L^2(S)}\big)^2
    +\big(\|H\|_{L^2(S)}^2-4M\big)\Big(1-\frac{2\|H_0\|_{L^2(S)}}{\|H\|_{L^2(S)}+2\sqrt M}\Big)+4\int_S (HH_0)_-\dd\calH^2\nonumber\\
    &\qquad\geq \big(2\sqrt M-|H_0||S|^{\frac{1}{2}}\big)^2
    +\big(\|H\|_{L^2(S)}^2-4M\big)\frac{2\sqrt M-|H_0||S|^{\frac{1}{2}}}{2\sqrt M}+4\int_S (HH_0)_-\dd\calH^2\,,
    \label{eq:PfMainMin-10}
\end{align}
where we have used \eqref{eq:PfMainMin-4} in the last inequality.
\end{step}

\begin{step}
    If $|H_0|\leq 2$ then \eqref{eq:PfMainMin-10} and $|S|\leq M$ yield
\begin{align}
    \int_S (H-H_0)^2\dd\calH^2
    &\geq (2-|H_0|)^2M
    +\big(\|H\|_{L^2(S)}^2-4M\big)\frac{2-|H_0|}{2}
    +4\int_S (HH_0)_-\dd\calH^2\,.
    \label{eq:PfMainMin-11}
\end{align}
Since $M\geq 4\pi$ we deduce for $0\leq H_0\leq 2$ that \eqref{eq:E(B)_minimal} holds.

Equality in \eqref{eq:E(B)_minimal} requires equality in \eqref{eq:PfMainMin-4} as well as $M=4\pi$, and therefore holds if and only if $S=\partial B$.
\end{step}

\begin{step}
    We finally consider the case $H_0<0$.
For $r\in (0,1]$, we have $\partial B_r(0)\subset \overline{B}$.
Sending $r\to 0$ in \eqref{eq:energy_sphere} we have
\begin{align}
    \calE_{H_0}(\partial B_r(0))\to 4\pi\,, \qquad \text{ as }r\to 0\,,
\end{align}
so that the infimum \eqref{eq:inf_ext} is not larger than $4\pi$.

On the other hand we have for any $S\in\calM$ with $S_+:=S\cap\{H\geq 0\}\cap\{K\geq 0\}$
\begin{align}
    \int_S (H-H_0)^2\dd\calH^2
    &\geq \int_{\{H\geq 0\}} (H-H_0)^2\dd\calH^2 \nonumber\\
    &\geq \int_{S_+} H^2\dd\calH^2
    + \int_{\{H\geq 0\}} (H_0^2+ 2|H_0|H)\dd\calH^2 > 16\pi\,,
    \label{eq:PfMainMin-6}
\end{align}
since $\int_{S_+} H^2 \dd\calH^2\geq 16\pi$ by an argument similar to the one used by Willmore \cite{Willmore_Riemannian}*{Lemma 7.2.1} in the proof that round spheres have minimal Willmore energy in the class of closed surfaces, see also \cite{RuppScharrer23}*{Remark 6.2}.
This proves the third claim.\qedhere
\end{step}
\end{proof}

We can relax the assumptions on $S$ and consider classes of generalized surfaces in the sense of oriented varifolds or Alexandrov immersed surfaces, see Section \ref{sec:OrVarif}.

\begin{proposition}\label{prop:rigidity_varifolds}
\leavevmode
\begin{enumerate}
\item Let $|H_0|\leq 2$.
Consider any oriented integral $2$-varifold $V=\underline{v}(S,\nu,\theta_+,\theta_-)$ with $\mu:=\|V\|$ supported in $\overline B$ and weak mean curvature $\vec H\in L^2(\mu)$.
Then for $M:=\max\{4\pi,|\Sigma|\}$ it holds
\begin{align}
    \calE_{H_0}(V) \geq
    &\frac{1}{4}(2-|H_0|)^2 M\,,
    \label{eq:LemmaHelfMLarge}
\end{align}
with equality if and only if $V=\underline{v}(\partial B,\sigma\nu_B,k,0)$ with $\sigma\in\{\pm 1\}$ such that $\sigma|H_0|=H_0$ and $k\in\N$.
\item Let $H_0<0$.
Then for any Alexandrov immersion $f\colon \Sigma \to \R^3$, $f(\Sigma)\subset \overline B$, for any varifold with enclosed volume \cite{RuppScharrer23}*{Hypothesis 4.5} and compact support in $\overline{B}$, and any volume varifold  \cite{scharrer2023properties}*{Definition 3.1} with compact support in $\overline{B}$, the value of $\calE_{H_0}$ is strictly larger than the infimum $4\pi$ in \eqref{eq:inf_ext}.
\end{enumerate}
\end{proposition}

\begin{proof}\leavevmode
\begin{enumerate}[wide]
\item Let $0<H_0\leq 2$.
As a consequence of \cite{MuRo14}*{Theorem 1}, we first remark that \eqref{eq:PfMainMin-4} and the characterization of the equality case also hold in the class of (oriented) varifolds with square-integrable weak mean curvature that are supported in $\overline B$, hence
\begin{equation}
    \frac{1}{4}\int_{\R^3} |\vec H|^2\dd\mu
    \geq M:=\max\{4\pi,\mu(\R^3)\}\,,
    \label{eq:PfMainMinGen-4}
\end{equation}
with equality for $\mu(\R^3)<4\pi$ if and only if $\mu=\calH^2\mres\partial B_r(x_0)$, $r\in (0,1]$, and for $|\Sigma|\geq 4\pi$ if and only if $\mu=k\calH^2\mres\partial B$ for some $k\in\N$.

We now represent $V$ in the form $V=\underline{v}(S,\nu^*,\theta^*_+,\theta^*_-)$ such that
\begin{equation}
    \vec H\cdot \nu^* = \sgn(H_0)|\vec H| \quad\mu\text{-a.e.}\,.
    \label{eq:CanRepIVo}
\end{equation}
In analogy to \eqref{eq:PfMainMin-11} we have
\begin{align}
    &4\calE_{H_0}(V)
    \geq (2-|H_0|)^2M
    +\big(\|H\|_{L^2(\mu)}^2-4M\big)\frac{2-|H_0|}{2}\sqrt M
    +4|H_0|\int_S |\vec{H}|\theta^*_-\dd\calH^2\,.
    \label{eq:PfMainMinGen-1}
\end{align}
Therefore
\begin{align}
    4\calE_{H_0}(V) \geq (2-|H_0|)^2M
\end{align}
with equality if and only if equality holds in \eqref{eq:PfMainMinGen-4} and $\theta^*_-=0$, hence if and only if $V=(\partial B,\nu_B,\theta_++\theta_-,0)$ with $\nu_B$ the inner normal of $B$.

\item Let $H_0<0$.
By considering shrinking spheres, the infimum in \eqref{eq:inf_ext} is at most $4\pi$.
In the class of varifolds with enclosed volume and compact support, we have $\calE_{H_0}>4\pi$ by \cite{RuppScharrer23}*{Corollary 4.11 (where there is a factor $1/(4\pi)$ missing for the first term on the right hand side)}.
For a volume varifold with compact support, the statement thus follows after employing \cite{scharrer2023properties}*{Proposition 3.6(6)(17).} For an Alexandrov immersed surface, the statement follows from \cite{RuppScharrer23}*{Theorem 1.5}.\qedhere
\end{enumerate}
\end{proof}

\subsection{The confinement problem for prescribed area}
We consider an additional area constraint and denote for given $a>0$ by $\calM_a$ the class of smooth, embedded closed surfaces $S\subset \overline{B}$ with surface area $\calH^2(S)=a$.
For $a>0$, $H_0\in\R$ we consider the minimization problem
\begin{align}
    \label{eq:inf_ext_area}
    \opt_a(H_0) := \inf\{ \calE_{H_0}(S) \mid S \in\calM_a\}\,.
\end{align}

In the case of small prescribed area and $H_0\in [0,2]$, minimizers are balls.
We retrieve rigidity of the ball for negative $H_0$, at least if Minkowski's inequality
\begin{align}\label{eq:minkowski}
    \int_S H\dd\calH^2 \geq \sqrt{16\pi |S|}
\end{align}
holds.
By \cite{MR2522433}*{Theorem 2} and \cite{MR3544938} this inequality is always satisfied for surfaces that are mean convex and star-shaped, or axi-convex.
\begin{proposition}\label{prop:rigidity_area_constraint}
Consider $a=4\pi r^2$ with $0<r\leq 1$.
\begin{enumerate}
    \item For $0\leq H_0\leq 2$ balls with radius $r$ are the unique minimizers in \eqref{eq:inf_ext_area}.
    \item Let $H_0<0$ and consider the class $\calM_{a,H_0}^*$ of surfaces $S\in\calM_a$ that satisfy in the case $H_0^2a>16\pi$ the Minkowski inequality
    \begin{equation}
        \int_S H\dd\calH^2 \geq 4\sqrt\pi\sqrt{a}\,,
        \label{eq:Minkowski-2}
    \end{equation}
    and in the case $H_0^2a\leq 16\pi$ a weak Minkowski inequality \begin{equation}
        \int_S H_+\dd\calH^2 \geq 4\sqrt\pi\sqrt{a}\,.
        \label{eq:WeakMinkowski}
    \end{equation}
    Then balls with radius $r$ are the unique minimizers of $\calE_{H_0}$ in $\calM_{a,H_0}^*$.
\end{enumerate}
\end{proposition}

\begin{proof}\leavevmode
\begin{enumerate}[wide]
\item In the case $0\leq H_0\leq 2$ we have by \eqref{eq:PfMainMin-10} (with $M=4\pi$) for any $S\in\calM_a$
\begin{align}
    \int_S (H-H_0)^2\dd\calH^2
    &\geq (2-H_0r)^24\pi
    +\big(\|H\|_{L^2(S)}^2-16\pi\big)\frac{2-H_0}{2}\sqrt{4\pi}
    +4H_0\int_S H_-\dd\calH^2\\
    &\geq \int_{\partial B_r} (H-H_0)^2\dd\calH^2\,,
\end{align}
which shows the optimality of balls.
Moreover, equality in this estimate implies equality in \eqref{eq:PfMainMin-4}, and hence that $S$ is a sphere of radius $r$.

\item Now let $H_0<0$ and first assume that \eqref{eq:Minkowski-2} holds.
Then we deduce similar as in \eqref{eq:PfMainMin-6} that
\begin{align}
    \int_S (H-H_0)^2\dd\calH^2
    &= \int_{S} H^2\dd\calH^2
    + \int_S (H_0^2 - 2 H_0H)\dd\calH^2 \nonumber\\
    &\geq 16\pi  + 4\pi H_0^2r^2 -16\pi H_0r = 4\pi(2-H_0r)^2
    = \int_{\partial B_r} (H-H_0)^2\dd\calH^2\,,
    \label{eq:PfMainMin-6a}
\end{align}
by \eqref{eq:Minkowski-2}.
This shows the optimality of balls.
Again, by the characterization of the equality case in Willmore's inequality, equality in this estimate implies that $S$ is a sphere of radius $r$.

If we have $H_0^2a\leq 16\pi$, then $r|H_0|\leq 2$, and we obtain from \eqref{eq:PfMainMin-10} and \eqref{eq:WeakMinkowski}
\begin{align}
    \int_S (H-H_0)^2\dd\calH^2
    &\geq (2-|H_0|r)^24\pi
    +4|H_0|\int_S H_+\dd\calH^2\\
    &\geq (2-|H_0|r)^24\pi
    +8|H_0|4\pi r\\
    &=(2-H_0r)^24\pi= \int_{\partial B_r} (H-H_0)^2\dd\calH^2\,. \qedhere
\end{align}
\end{enumerate}
\end{proof}

We finally consider prescribed area $a\geq 4\pi$ and $0\leq H_0< 2$.
We show that the Helfrich deficit strongly increases for prescribed area slightly above $4\pi$, see \cite{MuRo14} for a corresponding estimate for the Willmore energy.
This follows from an estimate for the (confined) Helfrich deficit by the (confined) Willmore deficit.
    
\begin{prop}\label{pro:est-impr-below}
Let $0\leq H_0<2$.
For all $a\geq 4\pi$ and all $S\in\calM_a$ it holds
\begin{align}
    \calE_{H_0}(S)-\frac{1}{4}(2-H_0)^2a
    &\geq \frac{2-H_0}{2}\big(\calE_{0}(S)-a\big)\,,
    \label{eq:CompDeficits2}\\
    \opt_a(H_0)-\opt_{4\pi}(H_0)
    &\geq \frac{2-H_0}{2}\big(\opt_a(0)-\opt_{4\pi}(0)\big)\,.
    \label{eq:CompDeficits2a}
\end{align}
Moreover, there exists a universal constant $c>0$ such that for all $a>4\pi$ sufficiently small (depending on $H_0$), we have
\begin{equation}
    \opt_a(H_0)-\opt_{4\pi}(H_0)
    \geq c\frac{2-H_0}{2}\sqrt{a-4\pi}\,.
    \label{eq:lb-4pi}
\end{equation}
\end{prop}

\begin{proof}
The inequality \eqref{eq:CompDeficits2} directly follows from \eqref{eq:PfMainMin-11}. As a consequence of \Cref{thm:main_ext_obstacle}, we have
\begin{align}
 \opt_{4\pi}(H_0) = (2-H_0)^2\pi\, .  \label{eq:opti_4pi} 
\end{align}
Passing to the infimum over all $S\in\calM_a$ on both sides of \eqref{eq:CompDeficits2} and using \eqref{eq:opti_4pi}
and $a\geq 4\pi$, we obtain \eqref{eq:CompDeficits2a}.

We further deduce from \eqref{eq:CompDeficits2}
\begin{align}
    &\int_S (H-H_0)^2\dd\calH^2 -(2-H_0)^24\pi\nonumber\\
    &\qquad \geq \big(\|H\|_{L^2(S)}^2-16\pi\big)\frac{2-H_0}{2}
    -(a-4\pi)2(2-H_0)+(2-H_0)^2(a-4\pi)\nonumber\\
    &\qquad =\big(\|H\|_{L^2(S)}^2-16\pi\big)\frac{2-H_0}{2}
    -H_0(2-H_0)(a-4\pi)\,.
    \label{eq:PfLB4Pi}
\end{align}
By \cite{MuRo14}*{Theorem 3} there exists $c'>0$ with $\opt_a(0)-4\pi\geq c'\sqrt{a-4\pi}$ for all $a>4\pi$ sufficiently small.
We therefore conclude with \eqref{eq:PfLB4Pi}
\begin{align}
    \int_S (H-H_0)^2\dd\calH^2 -(2-H_0)^24\pi
    &\geq 4 c'\frac{2-H_0}{2}\sqrt{a-4\pi}-H_0(2-H_0)(a-4\pi)\\
    \geq c\frac{2-H_0}{2}\sqrt{a-4\pi}
\end{align}
for all $a>4\pi$ sufficiently small.
Taking the infimum over all $S\in\calM_a$ implies \eqref{eq:lb-4pi}.
\end{proof}

\begin{remark}\label{rem:sharp_square_root}
The lower estimate \eqref{eq:lb-4pi} is sharp in its square root dependence on the area difference.

In fact one can review the construction in \cite{MuRo14}*{Proposition 2}.
There, for any $0<s<1$ sufficiently small, a surface $\Sigma_s$ inside $\overline{B}$ was constructed that only deviates from $\sphere^2$ in $B_s^2(0)\times\R^+$ and that has uniformly bounded mean curvature and small area difference to the unit sphere, more precisely
\begin{equation}
    |H_{\Sigma_s}| \leq C_1\,,\quad
    c_2s^4\leq |\Sigma_s|-4\pi \leq C_2s^4\qquad\text{ for all }0<s<s_0\,,
    \label{eq:ub}
\end{equation}
with $C_1,c_2,C_2>0$ independent of $s$ and $s_0>0$ sufficiently small.

For the difference in Helfrich energy, using the identity $\big|\Sigma_s\setminus\sphere^2\big| = |\Sigma_s|-4\pi+\big|\sphere^2\setminus\Sigma_s\big|$ in the second line, we obtain
\begin{align}
    \calE_{H_0}(\Sigma_s)-\calE_{H_0}(\sphere^2)
    &\leq \frac{1}{4}(C_1+H_0)^2\big|\Sigma_s\setminus\sphere^2\big|-\frac{1}{4}(2-H_0)^2\big|\sphere^2\setminus\Sigma_s\big|\\
    &=\frac{1}{4}(C_1+H_0)^2\big(|\Sigma_s|-4\pi\big)
    +\frac{1}{4}\Big((C_1+H_0)^2-\frac{1}{4}(2-H_0)^2\Big)\big|\sphere^2\setminus\Sigma_s\big|\\
    &\leq \frac{1}{4}(C_1+H_0)^2\big(|\Sigma_s|-4\pi\big) + C(C_1,H_0)s^2\leq C'(C_1,H_0,s_0)\sqrt{|\Sigma_s|-4\pi}
\end{align}
where we have used \eqref{eq:ub} in the first and the  last inequality.

Since $|\Sigma_s|>4\pi$ becomes arbitrarily close to $4\pi$ for $s>0$ sufficiently small, combined with \eqref{eq:opti_4pi} this proves an upper bound for $\opt_a(H_0)$ that matches the lower bound \eqref{eq:lb-4pi} with respect to the square root increase. 
\end{remark}

The following example shows that there are many examples of parameters $a,H_0$ such that the infimum $\opt_a(H_0)$ is zero.
The example also illustrates that it is necessary to pass to the set of immersed bubble trees when studying $\opt_a$, even when restricted to surfaces of spherical topology.
\begin{example}\label{ex:sphere_packing}
    For any $N\in\N, N\geq 2$, $\varepsilon>0$, there exists an embedded sphere $S=S_{N,\varepsilon}\subset B$ such that for $H_0\vcentcolon=2N$, we have $\H^2(S)= \frac{\pi}{4} N$ and $\calE_{H_0}(S)=o(1)$ as $\varepsilon\to 0$.
    In particular, $\opt_{N \pi/4}(2N)=0$.
\end{example}
\begin{proof}
    We may fit a cube $Q$ of side length $1/2$ inside $\overline{B}$.
    Divide $Q$ into $N^3$ smaller cubes of side length $1/(2N)$.
    Each of the smaller cubes contains a sphere of radius $r=1/(4N)$ such that each spheres touches its neighbors.
    Now, any two spheres may be connected by a small catenoidal neck of size $\varepsilon>0$ which only infinitesimally changes the area and the energy.
    It is possible to connect the lattice of spheres in such a way that the resulting surface $\Tilde{S}=\Tilde{S}_{N,\varepsilon}$ has spherical topology and satisfies
    \begin{align}
            \H^2(\Tilde{S})= \frac{N^3 \pi}{4N^2}+o(1)\,, \qquad
            \calE_{H_0}(\Tilde{S})=o(1)\,,\qquad\text{ as }\varepsilon\to 0\,.
        \end{align}
    Last, slightly deforming one of the spheres at the boundary of $Q$, we obtain a spherical surface $S$ which satisfies the area constraint $\H^2(S)=\frac{N^3 \pi}{4N^2}$.
\end{proof}

\begin{remark}\label{rem:sphere_packing}
    More generally, for $H_0>2$, any \emph{sphere packing} of $N$ spheres of equal radius $r=2/H_0<1$ in the unit ball is a trivial minimizer of $\calE_{H_0}$.
\end{remark}

\appendix

\section{Localized first variation of the Helfrich integrand}\label{sec:local_variaton}

\begin{lemma}\label{lem:local_first_vari}
    Let $\Phi\in \mathcal{F}_D$ be conformal, possibly with a single branch point at $z=0$. Let $\zeta\in L^\infty(D)$ be such that $\zeta\equiv 0$ in a neighborhood of $z=0$.
    Then for any $w\in C^\infty(D;\R^3)$ we have
    \begin{align}
         \left.\frac{\dd}{\dd t}\right\vert_{t=0} \calE(\Phi+tw; \zeta)  &= \int_D \frac{1}{2} (\vec H-H_0\nu)\cdot \Delta w\, \zeta\dd z \\
         &\qquad + \int \frac{1}{4}\left( (3 H - 4H_0)\nabla \nu - \vec H \times \nabla^\perp \nu -(2H_0 H -H_0^2)\nabla \Phi \right)\cdot \nabla w\,\zeta\dd z\,,
    \end{align} 
    and
    \begin{align}
         \left.\frac{\dd}{\dd t}\right\vert_{t=0} \mathcal{A}(\Phi+tw; \zeta) = \int_D \nabla \Phi\cdot\nabla w\, \zeta \dd z\,.
    \end{align}
\end{lemma}

\Cref{lem:local_first_vari} relies on the pointwise computations for the variation of mean curvature in \cite{MR2430975}, see also \cites{MR3524220,MR4076069}. The main difference is that we do not use integration by parts. For the convenience of the reader, we provide the details on how  \Cref{lem:local_first_vari} can be deduced from \cite{MR3524220}.

\begin{proof}[Proof of \Cref{lem:local_first_vari}]
    The set of $z\in D$ with $\zeta(z)\neq 0$ has positive distance from the branch point $z=0$ by assumption. Hence, we may rely on the pointwise computations in conformal coordinates from \cite{MR3524220}*{Proof of Theorem 4.57}. In particular, by \cite{MR3524220}*{(4.208)}, we have
    \begin{align}
        \left.\frac{\dd}{\dd t}\right\vert_{t=0} \dd \mu_{\Phi+tw} =  \nabla\Phi\cdot \nabla w\,\dd z\,.\label{eq:dt_mu}
    \end{align}
    For the mean curvature, combining \cite{MR3524220}*{(4.203),(4.206), and (4.207)}, with $\lambda$ the conformal factor of $\Phi$, we find
    \begin{align}
        &\left.\frac{\dd}{\dd t}\right\vert_{t=0}  H_{\Phi+tw} = - \left(\left.\frac{\dd}{\dd t}\right\vert_{t=0} g_{\Phi+tw}^{ij}\right) \nu \cdot \partial_{x^j} \Phi - g^{ij}\left[\partial_{x^i} \left(\left.\frac{\dd}{\dd t}\right\vert_{t=0} \nu_{\Phi+tw}\right)\cdot \partial_{x^j} \Phi + \partial_{x^i} \nu \cdot\partial_{x^j} w\right]\\
        &\quad  = e^{-4\lambda} \sum_{i,j} (\partial_{x^i} w\cdot \partial_{x^j}\Phi + \partial_{x^j} w\cdot\partial_{x^i}\Phi) \partial_{x^i} \nu \cdot \partial_{x^j} \Phi +e^{-2\lambda} \big(\partial_1 (\partial_1w\cdot\nu)+  e^{-2\lambda} \Delta w\cdot\nu \ \\
        &\quad = 2e^{-2\lambda} \nabla w\cdot\nabla \nu + e^{-2\lambda} \Delta w\cdot \nu\,,\label{eq:dt_H}
    \end{align}
    where we used the symmetry of the second fundamental form and the fact that $\partial_{x^i} \nu$ is tangent in the last step, cf.\ \cite{MR3524220}*{bottom of p.~367}. Combining \eqref{eq:dt_mu} and \eqref{eq:dt_H}, we obtain
    \begin{align}
        \left.\frac{\dd}{\dd t}\right\vert_{t=0} \Big( H_{\Phi+tw}\dd \mu_{\Phi+tw}\Big) = \Big(2 \nabla w\cdot \nabla \nu+ \Delta w\cdot\nu + H\,\nabla \Phi\cdot \nabla w\Big) \dd z\,.
    \end{align}
    For the localized Willmore energy, \eqref{eq:dt_mu} and \eqref{eq:dt_H} yield
    \begin{align}
        \left.\frac{\dd}{\dd t}\right\vert_{t=0} \Big(\frac{1}{4}(H_{\Phi+ tw})^2\dd \mu_{\Phi+tw}\Big) &= \frac{1}{2} H \left.\frac{\dd}{\dd t}\right\vert_{t=0} \Big( H_{\Phi+tw}\Big)e^{2\lambda}\dd z + \frac{1}{4}H^2 \nabla \Phi\cdot\nabla w\, \dd z \\
        &= \frac{1}{2} \vec H\cdot \Delta w \dd z+ \Big( H\nabla \nu + \frac{1}{4}H^2\nabla \Phi\Big)\cdot \nabla w\,\dd z \\
        & = \frac12 \vec H\cdot \Delta w \dd z + \Big(\frac{3}{4}H\nabla \nu - \frac{1}{4}\vec H\times \nabla^\perp\nu\Big)\cdot \nabla w\,\dd z\,, \label{eq:dt_H2}
    \end{align}
    where in the last step we used that by \cite{MR3524220}*{(1.60)} we have
    \begin{align}
        H^2 \nabla \Phi = - H\nabla \nu -\vec H\times \nabla^\perp \nu\,.
    \end{align}
    Combining \eqref{eq:dt_mu}, \eqref{eq:dt_H}, and \eqref{eq:dt_H2}, the statement follows.
\end{proof}

With the local statement at hand, we can now prove \Cref{lem:E_A_diff}.

\begin{proof}[Proof of \Cref{lem:E_A_diff}]
    We only consider the functional $\calE$, the argument for $\mathcal{A}$ is similar.  
    Let $|t|$ be sufficiently small, such that $\Phi+tw\in \mathcal{F}_{\Sigma}$. Let $U_i, \varphi_i$, and $\zeta_i$ be as in the definition of $\delta \E$ in \Cref{sec:prelims_variation}. Since $\delta\calE(\Phi)$ is independent of the particular charts used, we may assume that either $U_i\subset\subset \Sigma\setminus\{b_1, \dots, b_N\}$ or $w\vert_{U_i}\equiv 0$ for $i\in I$. Then
    \begin{align}
        \left.\frac{\dd}{\dd t}\right\vert_{t=0}\calE(\Phi+tw) = \sum_i  \left.\frac{\dd}{\dd t}\right\vert_{t=0}\calE(\Phi+tw;\zeta_i) = \sum_i \left.\frac{\dd}{\dd t}\right\vert_{t=0}\calE( (\Phi+tw)\circ\varphi_i;\zeta_i\circ\varphi_i)\,,
    \end{align}
    where in the second equality, we used that $\Phi+tw\vert_{U_i}$ is either constant in $t$ (if $U_i$ is a neighborhood of a branch point in which case $w\vert_{U_i}=0$) or  a Lipschitz immersion, in which case the functional is easily seen to be invariant with respect to composition with diffeomorphisms from the right (for instance, using that $W^{2,2}$-Lipschitz immersions can be strongly approximated by smooth immersions, see \cite{rupp2024global}*{Lemma 5.6} for a detailed proof). Now, using \Cref{lem:local_first_vari} and the definition of $\delta\calE$ in \eqref{eq:def_dE}, \eqref{eq:def_dE_chart}, and \eqref{eq:def_dE_partition_of_unity}, the formula \eqref{eq:E_A_diff_1} (for $\calE$) follows. 
\end{proof}

\section*{Acknowledgments}
This research was funded in whole, or in part, by the Austrian Science Fund (FWF), grant number \href{https://doi.org/10.55776/ESP557}{10.55776/ESP557}.
\bibliography{Lib}

\begin{bibdiv}
\begin{biblist}

\bib{MR458158}{article}{
      author={Adams, David~R.},
       title={A note on {R}iesz potentials},
        date={1975},
        ISSN={0012-7094,1547-7398},
     journal={Duke Math. J.},
      volume={42},
      number={4},
       pages={765\ndash 778},
         url={http://projecteuclid.org/euclid.dmj/1077311348},
      review={\MR{458158}},
}

\bib{MR143162}{article}{
      author={Alexandrov, Aleksandr~D.},
       title={A characteristic property of spheres},
        date={1962},
        ISSN={0003-4622},
     journal={Ann. Mat. Pura Appl. (4)},
      volume={58},
       pages={303\ndash 315},
         url={https://doi.org/10.1007/BF02413056},
      review={\MR{143162}},
}

\bib{Allard}{article}{
      author={Allard, William~K.},
       title={On the first variation of a varifold},
        date={1972},
        ISSN={0003-486X},
     journal={Ann. of Math. (2)},
      volume={95},
       pages={417\ndash 491},
         url={https://doi.org/10.2307/1970868},
      review={\MR{307015}},
}

\bib{MR3518329}{article}{
      author={Bernard, Yann},
       title={Noether's theorem and the {W}illmore functional},
        date={2016},
        ISSN={1864-8258,1864-8266},
     journal={Adv. Calc. Var.},
      volume={9},
      number={3},
       pages={217\ndash 234},
         url={https://doi.org/10.1515/acv-2014-0033},
      review={\MR{3518329}},
}

\bib{MR2843242}{article}{
      author={Bernard, Yann},
      author={Rivi\`ere, Tristan},
       title={Local {P}alais-{S}male sequences for the {W}illmore functional},
        date={2011},
        ISSN={1019-8385,1944-9992},
     journal={Comm. Anal. Geom.},
      volume={19},
      number={3},
       pages={563\ndash 599},
         url={https://doi.org/10.4310/CAG.2011.v19.n3.a5},
      review={\MR{2843242}},
}

\bib{MR3096502}{article}{
      author={Bernard, Yann},
      author={Rivi\`ere, Tristan},
       title={Singularity removability at branch points for {W}illmore
  surfaces},
        date={2013},
        ISSN={0030-8730,1945-5844},
     journal={Pacific J. Math.},
      volume={265},
      number={2},
       pages={257\ndash 311},
         url={https://doi.org/10.2140/pjm.2013.265.257},
      review={\MR{3096502}},
}

\bib{MR4706029}{article}{
      author={Bernard, Yann},
      author={Wheeler, Glen},
      author={Wheeler, Valentina-Mira},
       title={Analysis of the inhomogeneous {W}illmore equation},
        date={2024},
        ISSN={0294-1449,1873-1430},
     journal={Ann. Inst. H. Poincar\'e{} C Anal. Non Lin\'eaire},
      volume={41},
      number={1},
       pages={129\ndash 158},
         url={https://doi.org/10.4171/aihpc/79},
      review={\MR{4706029}},
}

\bib{blaschke}{book}{
      author={Blaschke, G.},
      author={Thomsen, G.},
       title={Vorlesungen {\"u}ber {D}ifferentialgeometrie und geometrische
  {G}rundlagen von {E}insteins {R}elativit{\"a}tstheorie {III}:
  {D}ifferentialgeometrie der {K}reise und {K}ugeln},
      series={Grundlehren der mathematischen Wissenschaften},
   publisher={Springer Berlin Heidelberg},
        date={2013},
        ISBN={9783642508233},
         url={https://books.google.de/books?id=0HmJBwAAQBAJ},
}

\bib{BoMM15}{article}{
      author={Bouzar, Lila},
      author={Menas, Ferhat},
      author={M\"uller, Martin~Michael},
       title={Toroidal membrane vesicles in spherical confinement},
        date={2015},
     journal={Phys. Rev. E},
      volume={92},
       pages={032721},
         url={https://link.aps.org/doi/10.1103/PhysRevE.92.032721},
}

\bib{Brakke}{book}{
      author={Brakke, Kenneth~A.},
       title={The motion of a surface by its mean curvature},
      series={Mathematical Notes},
   publisher={Princeton University Press, Princeton, N.J.},
        date={1978},
      volume={20},
        ISBN={0-691-08204-9},
      review={\MR{0485012}},
}

\bib{BrazdaLussardiStefanelli20}{article}{
      author={Brazda, Katharina},
      author={Lussardi, Luca},
      author={Stefanelli, Ulisse},
       title={Existence of varifold minimizers for the multiphase
  {C}anham-{H}elfrich functional},
        date={2020},
        ISSN={0944-2669,1432-0835},
     journal={Calc. Var. Partial Differential Equations},
      volume={59},
      number={3},
       pages={Paper No. 93, 26},
         url={https://doi.org/10.1007/s00526-020-01759-9},
      review={\MR{4098040}},
}

\bib{Canham}{article}{
      author={Canham, P.B.},
       title={The minimum energy of bending as a possible explanation of the
  biconcave shape of the human red blood cell},
        date={1970},
     journal={J. Theor. Biol.},
      volume={26},
      number={1},
       pages={61\ndash 81},
}

\bib{MR4076072}{article}{
      author={Da~Lio, Francesca},
      author={Palmurella, Francesco},
      author={Rivi\`ere, Tristan},
       title={A resolution of the {P}oisson problem for elastic plates},
        date={2020},
        ISSN={0003-9527,1432-0673},
     journal={Arch. Ration. Mech. Anal.},
      volume={236},
      number={3},
       pages={1593\ndash 1676},
         url={https://doi.org/10.1007/s00205-020-01499-2},
      review={\MR{4076072}},
}

\bib{DaDe18}{article}{
      author={Dall'Acqua, Anna},
      author={Deckelnick, Klaus},
       title={An obstacle problem for elastic graphs},
        date={2018},
        ISSN={0036-1410},
     journal={SIAM J. Math. Anal.},
      volume={50},
      number={1},
       pages={119\ndash 137},
}

\bib{MR3544938}{article}{
      author={Dalphin, J.},
      author={Henrot, A.},
      author={Masnou, S.},
      author={Takahashi, T.},
       title={On the minimization of total mean curvature},
        date={2016},
        ISSN={1050-6926,1559-002X},
     journal={J. Geom. Anal.},
      volume={26},
      number={4},
       pages={2729\ndash 2750},
         url={https://doi.org/10.1007/s12220-015-9646-y},
      review={\MR{3544938}},
}

\bib{DaMN18}{article}{
      author={Dayrens, Fran{\c{c}}ois},
      author={Masnou, Simon},
      author={Novaga, Matteo},
       title={Existence, regularity and structure of confined elasticae},
        date={2018},
        ISSN={1292-8119},
     journal={ESAIM, Control Optim. Calc. Var.},
      volume={24},
      number={1},
       pages={25\ndash 43},
}

\bib{DeulingHelfrich}{article}{
      author={Deuling, H.J.},
      author={Helfrich, W.},
       title={Red blood cell shapes as explained on the basis of curvature
  elasticity},
        date={1976},
     journal={Biophysical journal},
      volume={16},
      number={8},
       pages={861\ndash 868},
}

\bib{DoMR11}{article}{
      author={Dondl, Patrick~W.},
      author={Mugnai, Luca},
      author={Röger, Matthias},
       title={{Confined Elastic Curves}},
        date={2011},
        ISSN={0036-1399, 1095-712X},
     journal={SIAM J. Appl. Math.},
      volume={71},
      number={6},
       pages={2205\ndash 2226},
         url={https://doi.org/10.1137/100805339},
}

\bib{MR3896647}{article}{
      author={Eichmann, Sascha},
       title={The {H}elfrich boundary value problem},
        date={2019},
        ISSN={0944-2669,1432-0835},
     journal={Calc. Var. Partial Differential Equations},
      volume={58},
      number={1},
       pages={Paper No. 34, 26},
         url={https://doi.org/10.1007/s00526-018-1468-x},
      review={\MR{3896647}},
}

\bib{MR3099262}{book}{
      author={Giaquinta, Mariano},
      author={Martinazzi, Luca},
       title={An introduction to the regularity theory for elliptic systems,
  harmonic maps and minimal graphs},
     edition={Second},
      series={Appunti. Scuola Normale Superiore di Pisa (Nuova Serie) [Lecture
  Notes. Scuola Normale Superiore di Pisa (New Series)]},
   publisher={Edizioni della Normale, Pisa},
        date={2012},
      volume={11},
        ISBN={978-88-7642-442-7; 978-88-7642-443-4},
         url={https://doi.org/10.1007/978-88-7642-443-4},
      review={\MR{3099262}},
}

\bib{GrOk23}{article}{
      author={Grunau, Hans-Christoph},
      author={Okabe, Shinya},
       title={Willmore obstacle problems under {D}irichlet boundary
  conditions},
        date={2023},
        ISSN={0391-173X,2036-2145},
     journal={Ann. Sc. Norm. Super. Pisa Cl. Sci. (5)},
      volume={24},
      number={3},
       pages={1415\ndash 1462},
         url={https://doi.org/10.2422/2036-2145.202105_064},
      review={\MR{4675964}},
}

\bib{MR2522433}{article}{
      author={Guan, Pengfei},
      author={Li, Junfang},
       title={The quermassintegral inequalities for {$k$}-convex starshaped
  domains},
        date={2009},
        ISSN={0001-8708,1090-2082},
     journal={Adv. Math.},
      volume={221},
      number={5},
       pages={1725\ndash 1732},
         url={https://doi.org/10.1016/j.aim.2009.03.005},
      review={\MR{2522433}},
}

\bib{MR1913803}{book}{
      author={H\'elein, Fr\'ed\'eric},
       title={Harmonic maps, conservation laws and moving frames},
     edition={Second},
      series={Cambridge Tracts in Mathematics},
   publisher={Cambridge University Press, Cambridge},
        date={2002},
      volume={150},
        ISBN={0-521-81160-0},
         url={https://doi.org/10.1017/CBO9780511543036},
        note={Translated from the 1996 French original, With a foreword by
  James Eells},
      review={\MR{1913803}},
}

\bib{Helfrich}{article}{
      author={Helfrich, Wolfgang},
       title={Elastic properties of lipid bilayers: Theory and possible
  experiments},
        date={1973},
     journal={Zeitschrift für Naturforschung C},
      volume={28},
      number={11},
       pages={693–703},
}

\bib{MR1996773}{book}{
      author={H\"ormander, Lars},
       title={The analysis of linear partial differential operators. {I}},
      series={Classics in Mathematics},
   publisher={Springer-Verlag, Berlin},
        date={2003},
        ISBN={3-540-00662-1},
         url={https://doi.org/10.1007/978-3-642-61497-2},
      review={\MR{1996773}},
}

\bib{Hutc86}{article}{
      author={Hutchinson, John},
       title={{Second Fundamental Form for Varifolds and the Existence of
  Surfaces Minimising Curvature}},
        date={1986},
        ISSN={0022-2518},
     journal={Indiana U. Math. J.},
      volume={35},
      number={1},
       pages={45},
         url={https://doi.org/10.1512/iumj.1986.35.35003},
}

\bib{KaSM12a}{article}{
      author={Kahraman, Osman},
      author={Stoop, Norbert},
      author={M{\"u}ller, Martin~Michael},
       title={Fluid membrane vesicles in confinement},
        date={2012},
     journal={New Journal of Physics},
      volume={14},
      number={9},
       pages={095021},
}

\bib{KaSM12}{article}{
      author={Kahraman, Osman},
      author={Stoop, Norbert},
      author={M{\"u}ller, Martin~Michael},
       title={{Morphogenesis of membrane invaginations in spherical
  confinement}},
        date={2012},
        ISSN={0295-5075, 1286-4854},
     journal={Europhys. Lett.},
      volume={97},
      number={6},
       pages={68008},
         url={https://doi.org/10.1209/0295-5075/97/68008},
}

\bib{KellerMondinoRiviere2014}{article}{
      author={Keller, Laura Gioia~Andrea},
      author={Mondino, Andrea},
      author={Rivi\`ere, Tristan},
       title={Embedded surfaces of arbitrary genus minimizing the {W}illmore
  energy under isoperimetric constraint},
        date={2014},
        ISSN={0003-9527,1432-0673},
     journal={Arch. Ration. Mech. Anal.},
      volume={212},
      number={2},
       pages={645\ndash 682},
         url={https://doi.org/10.1007/s00205-013-0694-9},
      review={\MR{3176354}},
}

\bib{KubinLussardiMorandotti24}{article}{
      author={Kubin, Anna},
      author={Lussardi, Luca},
      author={Morandotti, Marco},
       title={Direct minimization of the {C}anham-{H}elfrich energy on
  generalized {G}auss graphs},
        date={2024},
        ISSN={1050-6926,1559-002X},
     journal={J. Geom. Anal.},
      volume={34},
      number={5},
       pages={Paper No. 121, 25},
         url={https://doi.org/10.1007/s12220-024-01564-2},
      review={\MR{4715387}},
}

\bib{KusnerMcGrath23}{article}{
      author={Kusner, Robert},
      author={McGrath, Peter},
       title={On the {C}anham problem: bending energy minimizers for any genus
  and isoperimetric ratio},
        date={2023},
        ISSN={0003-9527,1432-0673},
     journal={Arch. Ration. Mech. Anal.},
      volume={247},
      number={1},
       pages={Paper No. 10, 14},
         url={https://doi.org/10.1007/s00205-022-01833-w},
      review={\MR{4537344}},
}

\bib{MR2928715}{article}{
      author={Kuwert, Ernst},
      author={Li, Yuxiang},
       title={{$W^{2,2}$}-conformal immersions of a closed {R}iemann surface
  into {$\mathbb R^n$}},
        date={2012},
        ISSN={1019-8385,1944-9992},
     journal={Comm. Anal. Geom.},
      volume={20},
      number={2},
       pages={313\ndash 340},
         url={https://doi.org/10.4310/CAG.2012.v20.n2.a4},
      review={\MR{2928715}},
}

\bib{KuSc04}{article}{
      author={Kuwert, Ernst},
      author={Sch\"atzle, Reiner},
       title={{Removability of point singularities of Willmore surfaces}},
        date={2004},
        ISSN={0003-486X},
     journal={Ann. Math.},
      volume={160},
      number={1},
       pages={315\ndash 357},
         url={https://doi.org/10.4007/annals.2004.160.315},
      review={\MR{2119722 (2006c:58015)}},
}

\bib{MR2318282}{article}{
      author={Kuwert, Ernst},
      author={Sch\"atzle, Reiner},
       title={Branch points of {W}illmore surfaces},
        date={2007},
        ISSN={0012-7094,1547-7398},
     journal={Duke Math. J.},
      volume={138},
      number={2},
       pages={179\ndash 201},
         url={https://doi.org/10.1215/S0012-7094-07-13821-3},
      review={\MR{2318282}},
}

\bib{MR3059839}{article}{
      author={Kuwert, Ernst},
      author={Sch\"atzle, Reiner},
       title={Closed surfaces with bounds on their {W}illmore energy},
        date={2012},
        ISSN={0391-173X,2036-2145},
     journal={Ann. Sc. Norm. Super. Pisa Cl. Sci. (5)},
      volume={11},
      number={3},
       pages={605\ndash 634},
      review={\MR{3059839}},
}

\bib{MR3156649}{book}{
      author={Marcus, Moshe},
      author={V\'eron, Laurent},
       title={Nonlinear second order elliptic equations involving measures},
      series={De Gruyter Series in Nonlinear Analysis and Applications},
   publisher={De Gruyter, Berlin},
        date={2014},
      volume={21},
        ISBN={978-3-11-030515-9; 978-3-11-030531-9},
      review={\MR{3156649}},
}

\bib{MR2989995}{article}{
      author={Mondino, Andrea},
      author={Rivi\`ere, Tristan},
       title={Willmore spheres in compact {R}iemannian manifolds},
        date={2013},
        ISSN={0001-8708,1090-2082},
     journal={Adv. Math.},
      volume={232},
       pages={608\ndash 676},
         url={https://doi.org/10.1016/j.aim.2012.09.014},
      review={\MR{2989995}},
}

\bib{MR3276119}{article}{
      author={Mondino, Andrea},
      author={Rivi\`ere, Tristan},
       title={Immersed spheres of finite total curvature into manifolds},
        date={2014},
        ISSN={1864-8258,1864-8266},
     journal={Adv. Calc. Var.},
      volume={7},
      number={4},
       pages={493\ndash 538},
         url={https://doi.org/10.1515/acv-2013-0106},
      review={\MR{3276119}},
}

\bib{MR4076069}{article}{
      author={Mondino, Andrea},
      author={Scharrer, Christian},
       title={Existence and regularity of spheres minimising the
  {C}anham-{H}elfrich energy},
        date={2020},
        ISSN={0003-9527,1432-0673},
     journal={Arch. Ration. Mech. Anal.},
      volume={236},
      number={3},
       pages={1455\ndash 1485},
         url={https://doi.org/10.1007/s00205-020-01497-4},
      review={\MR{4076069}},
}

\bib{MondinoScharrer2023}{article}{
      author={Mondino, Andrea},
      author={Scharrer, Christian},
       title={A strict inequality for the minimization of the {W}illmore
  functional under isoperimetric constraint},
        date={2023},
        ISSN={1864-8258,1864-8266},
     journal={Adv. Calc. Var.},
      volume={16},
      number={3},
       pages={529\ndash 540},
         url={https://doi.org/10.1515/acv-2021-0002},
      review={\MR{4609797}},
}

\bib{Muel19}{article}{
      author={M{\"u}ller, Marius},
       title={An obstacle problem for elastic curves: existence results},
        date={2019},
        ISSN={1463-9963},
     journal={Interfaces Free Bound.},
      volume={21},
      number={1},
       pages={87\ndash 129},
}

\bib{MR1366547}{article}{
      author={M\"uller, Stefan},
      author={\v{S}ver\'ak, Vladim\'ir},
       title={On surfaces of finite total curvature},
        date={1995},
        ISSN={0022-040X,1945-743X},
     journal={J. Differential Geom.},
      volume={42},
      number={2},
       pages={229\ndash 258},
         url={http://projecteuclid.org/euclid.jdg/1214457233},
      review={\MR{1366547}},
}

\bib{MuRo14}{article}{
      author={Müller, Stefan},
      author={Röger, Matthias},
       title={{Confined structures of least bending energy}},
        date={2014},
        ISSN={0022-040X},
     journal={J. Differ. Geom.},
      volume={97},
      number={1},
       pages={109\ndash 139},
         url={https://doi.org/10.4310/jdg/1404912105},
      review={\MR{3229052}},
}

\bib{Pozz23}{article}{
      author={Pozzetta, Marco},
       title={Confined {W}illmore energy and the area functional},
        date={2023},
        ISSN={1019-8385,1944-9992},
     journal={Comm. Anal. Geom.},
      volume={31},
      number={2},
       pages={407\ndash 447},
         url={https://doi.org/10.4310/cag.2023.v31.n2.a7},
      review={\MR{4685027}},
}

\bib{RiPK14}{article}{
      author={Rim, Jee~E.},
      author={Purohit, Prashant~K.},
      author={Klug, William~S.},
       title={Mechanical collapse of confined fluid membrane vesicles},
        date={2014},
        ISSN={1617-7940},
     journal={Biomechanics and Modeling in Mechanobiology},
      volume={13},
      number={6},
       pages={1277\ndash 1288},
         url={https://doi.org/10.1007/s10237-014-0572-x},
}

\bib{MR2430975}{article}{
      author={Rivi\`ere, Tristan},
       title={Analysis aspects of {W}illmore surfaces},
        date={2008},
        ISSN={0020-9910,1432-1297},
     journal={Invent. Math.},
      volume={174},
      number={1},
       pages={1\ndash 45},
         url={https://doi.org/10.1007/s00222-008-0129-7},
      review={\MR{2430975}},
}

\bib{MR3276154}{article}{
      author={Rivi\`ere, Tristan},
       title={Variational principles for immersed surfaces with {$L^2$}-bounded
  second fundamental form},
        date={2014},
        ISSN={0075-4102,1435-5345},
     journal={J. Reine Angew. Math.},
      volume={695},
       pages={41\ndash 98},
         url={https://doi.org/10.1515/crelle-2012-0106},
      review={\MR{3276154}},
}

\bib{MR3524220}{incollection}{
      author={Rivi\`ere, Tristan},
       title={Weak immersions of surfaces with {$L^2$}-bounded second
  fundamental form},
        date={2016},
   booktitle={Geometric analysis},
      series={IAS/Park City Math. Ser.},
      volume={22},
   publisher={Amer. Math. Soc., Providence, RI},
       pages={303\ndash 384},
         url={https://doi.org/10.1090/pcms/022/07},
      review={\MR{3524220}},
}

\bib{RuppScharrer23}{article}{
      author={Rupp, Fabian},
      author={Scharrer, Christian},
       title={Li-{Y}au inequalities for the {H}elfrich functional and
  applications},
        date={2023},
        ISSN={0944-2669},
     journal={Calc. Var. Partial Differential Equations},
      volume={62},
      number={2},
       pages={Paper No. 45},
         url={https://doi.org/10.1007/s00526-022-02381-7},
      review={\MR{4525726}},
}

\bib{rupp2024global}{misc}{
      author={Rupp, Fabian},
      author={Scharrer, Christian},
       title={Global regularity of integral 2-varifolds with square integrable
  mean curvature},
        date={2024},
        note={To appear in J. Math. Pures Appl.},
}

\bib{Scharrer22NLA}{article}{
      author={Scharrer, Christian},
       title={Embedded {D}elaunay tori and their {W}illmore energy},
        date={2022},
        ISSN={0362-546X,1873-5215},
     journal={Nonlinear Anal.},
      volume={223},
       pages={Paper No. 113010, 23},
         url={https://doi.org/10.1016/j.na.2022.113010},
      review={\MR{4438233}},
}

\bib{scharrer2023properties}{misc}{
      author={Scharrer, Christian},
       title={Properties of surfaces with spontaneous curvature},
        date={2023},
}

\bib{scharrer2025energyquantizationconstrainedwillmore}{misc}{
      author={Scharrer, Christian},
      author={West, Alexander},
       title={Energy quantization for constrained {W}illmore surfaces},
        date={2025},
         url={https://arxiv.org/abs/2505.19898},
}

\bib{MR4865096}{article}{
      author={Scharrer, Christian},
      author={West, Alexander},
       title={On the minimization of the {W}illmore energy under a constraint
  on total mean curvature and area},
        date={2025},
        ISSN={0003-9527,1432-0673},
     journal={Arch. Ration. Mech. Anal.},
      volume={249},
      number={2},
       pages={Paper No. 17, 62},
         url={https://doi.org/10.1007/s00205-025-02087-y},
      review={\MR{4865096}},
}

\bib{MR2472179}{article}{
      author={Sch\"{a}tzle, Reiner},
       title={Lower semicontinuity of the {W}illmore functional for currents},
        date={2009},
        ISSN={0022-040X,1945-743X},
     journal={J. Differential Geom.},
      volume={81},
      number={2},
       pages={437\ndash 456},
         url={http://projecteuclid.org/euclid.jdg/1231856266},
      review={\MR{2472179}},
}

\bib{MR2928137}{article}{
      author={Schygulla, Johannes},
       title={Willmore minimizers with prescribed isoperimetric ratio},
        date={2012},
        ISSN={0003-9527,1432-0673},
     journal={Arch. Ration. Mech. Anal.},
      volume={203},
      number={3},
       pages={901\ndash 941},
         url={https://doi.org/10.1007/s00205-011-0465-4},
      review={\MR{2928137}},
}

\bib{Simon}{article}{
      author={Simon, Leon},
       title={Existence of surfaces minimizing the {W}illmore functional},
        date={1993},
        ISSN={1019-8385},
     journal={Comm. Anal. Geom.},
      volume={1},
      number={2},
       pages={281\ndash 326},
         url={https://doi.org/10.4310/CAG.1993.v1.n2.a4},
      review={\MR{1243525}},
}

\bib{MR1650335}{article}{
      author={Topping, Peter},
       title={Mean curvature flow and geometric inequalities},
        date={1998},
        ISSN={0075-4102,1435-5345},
     journal={J. Reine Angew. Math.},
      volume={503},
       pages={47\ndash 61},
         url={https://doi.org/10.1515/crll.1998.099},
      review={\MR{1650335}},
}

\bib{Willmore65}{article}{
      author={Willmore, Thomas~J.},
       title={Note on embedded surfaces},
        date={1965},
     journal={An. \c{S}ti. Univ. ``Al. I. Cuza" Ia\c{s}i Sec\c{t}. I a Mat.
  (N.S.)},
      volume={11B},
       pages={493\ndash 496},
}

\bib{Willmore_Riemannian}{book}{
      author={Willmore, Thomas~J.},
       title={Riemannian geometry},
      series={Oxford Science Publications},
   publisher={The Clarendon Press, Oxford University Press, New York},
        date={1993},
        ISBN={0-19-853253-9},
      review={\MR{1261641}},
}

\bib{Wojt17}{article}{
      author={Wojtowytsch, Stephan},
       title={Helfrich's energy and constrained minimisation},
        date={2017},
        ISSN={1539-6746,1945-0796},
     journal={Commun. Math. Sci.},
      volume={15},
      number={8},
       pages={2373\ndash 2386},
         url={https://doi.org/10.4310/CMS.2017.v15.n8.a10},
      review={\MR{3744688}},
}

\bib{Wojt21}{article}{
      author={Wojtowytsch, Stephan},
       title={Confined elasticae and the buckling of cylindrical shells},
        date={2021},
        ISSN={1864-8258},
     journal={Adv. Calc. Var.},
      volume={14},
      number={4},
       pages={555\ndash 587},
}

\end{biblist}
\end{bibdiv}
\bibliographystyle{abbrev}
\end{document}